\documentclass[review,3p,10pt]{elsarticle}
\biboptions{sort&compress}

\usepackage{graphicx}
\usepackage{epstopdf}
\usepackage{amsmath}
\usepackage{cases}
\usepackage[caption=false,font=footnotesize]{subfig}
\usepackage{fixltx2e}
\usepackage{amssymb}
\usepackage{mathtools}
\usepackage{bm}
\usepackage{multirow}
\usepackage{color}
\usepackage{xcolor}
\usepackage{algorithm}
\usepackage{algorithmic}
\usepackage{amsthm}
\usepackage{mathrsfs}
\usepackage{textcomp,mathcomp}
\usepackage{threeparttable}
\usepackage{bbding}
\usepackage{pifont}
\usepackage{multicol}
\usepackage{enumitem}
\usepackage{adjustbox}

\usepackage{makecell}

\newtheorem{remark}{Remark}
\usepackage{lineno}

\usepackage{adjustbox}

\usepackage{url}

\newcommand{\tabincell}[2]{\renewcommand\arraystretch{0.8}\begin{tabular}{@{}#1@{}}#2\end{tabular}}

\journal{International Journal of Hydrogen Energy}

\begin{document}
\begin{frontmatter}

\title{Dynamic Modeling and Control of Multi-Stack Alkaline Water Electrolysis Systems with Shared Gas Separators and Lye Circulation:\\
{\color{black}Industrial Data-Based Validation and Simulation }
}


\author[label1]{Yiwei~Qiu}
\author[label1]{Jiatong~Li}
\author[label1]{Yangjun~Zeng}
\author[label1]{Yi~Zhou\corref{cor1}}
\ead{yizhou3230@scu.edu.cn}
\author[label1]{Shi~Chen}
\author[label1]{Xiaoyan~Qiu}
\author[label1]{Buxiang~Zhou}
\author[label2]{Ge~He}
\author[label2]{Xu~Ji}
\author[label3]{Wenying~Li}

\address[label1]{College of Electrical Engineering, Sichuan University, Chengdu, 610065, China}
\cortext[cor1]{Corresponding author}
\address[label2]{School of Chemical Engineering, Sichuan University, Chengdu, 610065, China}
\address[label3]{Tsinghua Sichuan Energy Internet Research Institute, Chengdu, 610213, China}

\begin{abstract}
  An emerging approach for large-scale renewable hydrogen production is integrating multiple alkaline water electrolysis (AWE) stacks into one balance-of-plant (BoP) system, sharing gas-lye separation and lye circulation components. While this configuration, termed $N$-in-1, reduces cost and complexity, its dynamic performance under fluctuating power remains unclear compared with conventional 1-in-1 systems. This paper develops a state-space model of the multi-stack AWE system, capturing lye circulation, temperature, and hydrogen-to-oxygen (HTO) dynamics, calibrated via experiments on a 4,000 Nm$^3$/h-rated 4-in-1 system. A {\color{black} mixed-integer quadratic programming (MIQP)-based predictive controller} is then designed to coordinate inter-stack current distribution, lye flow, and cooling for load tracking and operational stability. {\color{black}Simulations on the experimentally validated model} show that a $4$-in-1 system achieves similar performance compared to four parallel 1-in-1 systems {\color{black}under continuous operation}. Differences in load-tracking, temperature stabilization errors, and specific energy consumption remain below 0.015 MW, 0.346 K, and 0.001 kWh/Nm$^3$ under wind power supply {\color{black}when all stacks remain online}.

\end{abstract}


\end{frontmatter}

{\color{black}
\section*{Nomenclature}

\addcontentsline{toc}{section}{Nomenclature}
\begin{multicols}{2}
	
	\footnotesize
	\setlength{\columnsep}{18pt}

	\noindent\textbf{Abbreviations}
	\begin{description}[leftmargin=2.2cm, labelwidth=2.0cm, labelsep=0.2cm, align=left, itemsep=0pt, topsep=0pt, font=\normalfont]
		\item[ReP2H] Renewable power-to-hydrogen
		\item[AWE] Alkaline water electrolysis
		\item[BoP] Balance of plant
		\item[HTO] Hydrogen-to-oxygen
		\item[OTH] Oxygen-to-hydrogen
		\item[MPC]	Model predictive control
		\item[NMPC] Nonlinear MPC
		\item[SMPC] Stochastic MPC
		\item[MIQP] Mixed-integer quadratic programming
		\item[MIMO]	Multi-input multi-output
		\item[MPPT] Maximal power point tracking
		\item[AGC] Automatic generation control
	   	\item[DD]  Double description
		\item[RMSE] Root-mean-square error
	\end{description}
	\vspace{5pt}

	\noindent\textbf{Indices}
	
	\begin{description}[leftmargin=2.3cm, labelwidth=2.1cm, labelsep=0.2cm, align=left, itemsep=0pt, topsep=0pt, font=\normalfont]
		\addcontentsline{toc}{section}{Nomenclature}
		\item[$i, j$] Index for stacks/pumps/components
		\item[$k$] Index for time intervals
	\end{description}
	\vspace{5pt}

	\noindent\textbf{Variables}
			
	\begin{description}[leftmargin=2.3cm, labelwidth=2.1cm, labelsep=0.2cm, align=left, itemsep=0pt, topsep=0pt, font=\normalfont]
		\addcontentsline{toc}{section}{Nomenclature}
		\item[$U_i^{{\rm{cell}}}$, $I_i$] Cell voltage and current of the $i$th stack
		\item[$P_i$] Power consumption of the $i$th stack
		\item[$T_{\text{s},i}$] Temperature of the $i$th stack
		\item[$\eta_i^{\rm{cell}}$] Faraday efficiency of the $i$th stack
		\item[$\dot{n}_{i}^{\text{H}_2,\text{prod}}$]  Hydrogen production rate of the $i$th stack
		\item[$\dot{n}_{i}^{\text{O}_2,\text{prod}}$] Oxygen production rate of the $i$th stack
		
		\item[$\Delta p_i$]  Pressure drop across the $i$th stack
		\item[$R_i^{\text{ca/an}}$]  Cathode/anode-side flow resistance
		\item[$v_{i,\text{gas/lye/mix}}^{\text{ca/an}}$] Gas/lye/mixture flow rate in the $i$th stack
		\item[$v_{\text{tot,lye}}$] Total lye flow through the BoP
		\item[$V_{i,\text{gas/lye}}^{\text{ca/an}}$]  Volume of gas/lye in the $i$th stack
		\item[$\alpha_{i}^{\text{ca/an}}$] Cathode/anode-side volumetric gas ratio
				
		\item[$T_{i,\text{s,in/out}}$] Inlet/outlet lye temperature of the $i$th stack
		\item[$T_{\text{sep,out}}$] Separator outlet temperature
		\item[$T_{\text{c,out}}$] Cooling water outlet temperature
		\item[$\Delta T$] Logarithmic temperature difference in the heat exchanger
		\item[$C_{i,\text{s}}$] Total heat capacity of the $i$th stack and lye
		\item[$C_{i,\text{lye,s}}$] Heat capacity of lye in the $i$th stack
		\item[$V_{i,\text{s,lye}}$] Volume of lye in the $i$th stack
		\item[$Q_{i,\text{ele}}$] Stack heating flow
		\item[$Q_{i,\text{s,diss}}$] Heat dissipation from the $i$th stack
		\item[$Q_{i,\text{s,conv}}$] Convective heat loss from the $i$th stack
		\item[$Q_{i,\text{s,rad}}$] Radiative heat loss from the $i$th stack
		\item[$Q_{\text{sep,conv/rad}}$] Separator convective/radiative heat loss
		\item[$v_{\text{c}}$]  Cooling water flow rate
		
		\item[HTO] Gas-phase HTO impurity concentration
		
		\item[$\dot{n}^{\text{H}_2,\text{im}}_i$] HTO impurity crossover from the cathode side to anode side in the $i$th stack
		
		\item[$\dot{n}^{\rm{H_2},\rm{lye/diff/conv}}$] Molar flow of HTO impurity brought in by lye circulation/diffusion/convection
		\item[$n^{\text{H}_2,\text{an}}_{i}$] Anode half-cell HTO molar quantity
		\item[$n^{\text{H}_2,\text{sep}}_{\text{liq/gas}}$] Separator liquid/gas-phase molar quantity of HTO impurity
		\item[$\dot{n}^{\text{H}_2,\text{im}}_{i,1}$]  HTO molar flow from the $i$th stack to the separator
		\item[$\dot{n}^{\text{H}_2,\text{im}}_2$]  HTO molar flow to the separator gas phase
		\item[$\dot{n}^{\text{H}_2,\text{im}}_{\text{out}}$]  System outlet HTO molar flow
		\item[$V_{\rm{lye}}^{\rm{an}}$] Lye volume in the anode-side half-cell

		\item[$\beta^{I}_{i,k}$] The $k$th binary variable used for discretizing $I$
		\item[$\delta^{I,T}_{i,k}$] Intermediate variable in mixed-integer reformulation

	\end{description}
	\vspace{5pt}

	\noindent\textbf{Parameters}
	
	\begin{description}[leftmargin=2.3cm, labelwidth=2.1cm, labelsep=0.2cm, align=left, itemsep=0pt, topsep=0pt, font=\normalfont]
		\addcontentsline{toc}{section}{Nomenclature}
		\item[$N^\text{cell}$] Number of cells in a stack
		\item[$\phi_{s}$] Diameter of the stack
		\item[$l$, $d$] Equivalent length and diameter of stack flow path
			
		\item[$U_i^{\rm{rev}}$, $U^\text{th}$]  Reversible and thermoneutral voltages
		\item[$r_{i,1}$, $r_{i,2}$, $r_{i,3}$, $s_{i}$]			
		\item[$t_{i,1}$, $t_{i,2}$, $t_{i,3}$]  Constant parameters in cell voltage model
		\item[$f_{i,1}$, $f_{i,2}$]  Constant parameters in the Faraday efficiency model
		\item[$\rho$] System operating pressure
		\item[$\Delta \rho$]  Cathode-anode pressure difference		
	
		\item[$\mu_{\text{lye}}$, $\mu_{\text{H}_2/\text{O}_2}$] Dynamic viscosity of lye and gas
		
		\item[$\rho_\text{lye}$] Density of lye
		\item[$\rho_{\text{c}}$] Density of cooling water
		
		\item[$c_\text{lye}$] Specific heat capacity of lye
		\item[$c_{\text{c}}$] Specific heat capacity of cooling water
		
		\item[$T_\text{am}$] Ambient temperature
		\item[$T_{\text{c,in}}$] Inlet temperature of cooling water
		
		\item[$C_{\text{sep}}$] Total heat capacity of separator
		\item[$C_{\text{he}}$] Total heat capacity of heat exchanger
		\item[$C_{\text{c}}$] Total heat capacity of cooling coil
		
		\item[$A_{\text{s,diss}}$] Heat dissipation area of the stack
		\item[$A_\text{c}$] Heat exchange area of heat exchanger
		\item[$k$] Heat transfer coefficient of heat exchanger

		\item[$F$] The Faraday's constant		
		\item[$R$] Gas constant
		\item[$\sigma$] The Stefan-Boltzmann constant
		
		\item[$\epsilon_\text{s}$] Emissivity of the stack
		\item[$h_\text{s}$] Stack natural convection coefficient
		
		\item[$S^{\rm{H_2}}_\text{eff}$] Effective hydrogen carrying coefficient in the  lye
		
		\item[$D^{\rm{H_2}}_{\rm{eff}}$]  Effective diffusion coefficient
		\item[$\Delta c^{\rm{H_2}}$]  Differential hydrogen concentration
		\item[$\delta$]  Diaphragm thickness
		\item[$K^{\text{H}_2}_{\text{eff}}$] Effective hydrogen permeability of the diaphragm

		\item[$\tau_{\rm{sep}}$] Separation time constant
		\item[$V^{\text{sep}}_{\text{gas}}$] Separator gas phase volume

		\item[$P_{\text{tot}}^{\text{ref}}$] Power consumption reference
		\item[$T_{\text{s,out}}^{\text{ref}}$] Temperature reference for the $i$th stack
		
		\item[$\lambda^{\text{prod}}$, $\lambda^{\text{track}}$, $\lambda^{\text{temp}}$]
		\item[$\lambda^{\text{I}}$, $\lambda^{\text{lye}}$, $\lambda^{\text{c}}$]  Weight factors in control objective

		\item[$\overline{T}$] Stack temperature limit
		\item[$\overline{I}$, $\overline{P}^\text{ele}_i$] Rectifier current and power limits
		\item[$\overline{U}^\text{cell}$] Cell voltage limit
		\item[$r^{\text{H}_2,\text{prod,up/down}}$] Load ramping limits
		
		\item[$\underline{v}_{\text{lye}}$, $\overline{v}_{\text{lye}}$] Lye flow limits
		\item[$\underline{v}_{\text{c}}$, $\overline{v}_{\text{c}}$] Coolant flow limits
		\item[$\bm{A}$, $\bm{b}$] Constant matrix and vector in DD representation of production model
	
		\item[$N^\text{d}$] Number of binary bits in discretization
		\item[$\Delta I$] Discretization step size of stack current
		\item[$M$] A sufficiently large constant for mixed-integer reformulation
		
	\end{description}
	
\end{multicols}
}

\section{Introduction}
\label{sec:intro}


Water electrolysis for hydrogen production is a key pathway for renewable energy utilization and low-carbon transition of energy, transport, and chemical sectors, with rapid global deployment \cite{hosseini2022hydrogen}. Many countries have incorporated renewable power-to-hydrogen (ReP2H) into their strategies, targeting 520 GW of installed capacity by 2030 \cite{IEA2024Global}. Among available technologies, alkaline water electrolysis (AWE) stands out for its maturity, large capacity, long lifespan, and low cost, making it a preferred choice in many ReP2H projects \cite{emam2024review,zeng2024investment}.

Although the capacity of a single AWE unit has increased from 1,000 Nm$^3$/h (5 MW) in 2020 to 3,000 Nm$^3$/h (15 MW) \cite{emam2024review}, it remains limited for large-scale deployment. Hydrogen plants therefore require many electrolyzers, resulting in large land use, complex system layouts, and high investment costs \cite{emam2024review}. The costs of electrolyzers, land, and infrastructure are critical in practice \cite{lange2024modularization,niblett2024review}. To improve scalability and reduce costs, shared balance-of-plant (BoP) designs have emerged \cite{lange2024modularization,niblett2024review,chi2024elevating,zhang2024optimizing,yang2022application,rizwan2021design,chen2024enhancing,shi2023plant}.

A conventional AWE system consists of one electrolysis stack and a BoP system, including lye-gas separation, cooling, lye circulation, and auxiliary components such as pumps, valves, and control systems, usually integrated as a skid-mounted unit \cite{emam2024review,niblett2024review}. This configuration is referred to as a \emph{1-in-1} system. In large plants, however, the 1-in-1 design leads to excessive land use, complex topology, and high cost. To address this, \emph{$N$-in-1} configurations (also termed \emph{multi-electrolyzer} systems \cite{liang2024large,li2023exploration}) allow multiple stacks to share a single BoP system, reducing system complexity and cost \cite{niblett2024review}.

Industrial practice has demonstrated that replacing a 1-in-1 configuration with a 2-in-1 system reduces costs by 25\% and land and utility requirements by 45\%, while a 4-in-1 system achieves reductions of 35\% and 55\%, respectively \cite{yang2022application}. Consequently, many large ReP2H projects adopt $N$-in-1 configurations. For example, the \emph{Baofeng Energy Green Hydrogen Development Project} and \emph{Songwon Hydrogen Energy Industrial Park Project} use 2-in-1 systems, while the \emph{Sinopec Kuqa Green Hydrogen Demonstration Project} and the \emph{Da'an Wind and Solar Green Hydrogen Synthesis Ammonia Integration Demonstration Project} adopt 4-in-1 configurations \cite{meng2024advantages,zhangalkaline,zhai2024review}.
These developments raise an important question regarding the operational performance of $N$-in-1 systems under fluctuating renewable power inputs.

With the growth of ReP2H projects in recent years, research on AWE systems has expanded rapidly.
At the component level, many studies focus on electrolysis stacks, including electrochemical performance \cite{sanchez2018semi,huang2022multiphysics}, degradation \cite{esfandiari2024metal}, flow fields \cite{wang2023non}, and spatial distributions of temperature \cite{zhang2024temperature}, bubbles \cite{cao2024investigation}, and stray currents \cite{qi2023channel,sakas2024influence,dogan2024experimental}. Other work addresses BoP components, such as separator design and efficiency analysis \cite{wang2022experimental,zeng2020experimental,hu2024study,li2026structural}, heat exchangers \cite{rashidi2022progress,jiang2024performance}, pressure loss and energy consumption in lye circulation  \cite{qi2023channel}, and compression and purification systems \cite{kang2018optimal}, and Luo et al. \cite{luo2026hierarchical} present a data-driven framework to identify the causal interactions among BoP components.

Based on these studies, system-level models and controllers have been developed. Ulleberg et al. \cite{ulleberg2003modeling} proposed a lumped dynamic model of electrochemical and thermal behavior. David et al. \cite{david2020dynamic} developed a high-order mass transfer model including pressure, liquid level, and composition dynamics. Qi et al. \cite{qi2023thermal,qi2023design,qi2021pressure} introduced state-space models for temperature dynamics and HTO impurity accumulation under variable loads, along with thermal and pressure management strategies. Chen et al. \cite{cheng2024self} proposed an adaptive MPPT strategy considering temperature and pulse rectification effects. Li et al. \cite{li2022active} developed a control framework integrating temperature, pressure, and lye circulation. These models have been validated experimentally at scales from 500 W to 3 MW \cite{li2022study,ding2023experimental,sakas2022dynamic}. Qiu et al. \cite{qiu2023dynamic} further proposed online parameter identification for real-time calibration and diagnostics. Reviews of AWE modeling and control can be found in \cite{olivier2017low,daoudi2024overview,aguirre2024control}.

\begin{table}[tb]\scriptsize
	\renewcommand{\arraystretch}{1.22}
	\caption{{\color{black}Summary of the latest research related to process modeling and control of AWE systems}}\vspace{6pt}
	\label{tab:literature}
	\centering
	\begin{adjustbox}{width=\textwidth}
		\begin{tabular}{cccccccccc}
			\hline \hline
			\multirow{2}{*}{Literature}                   & \multirow{2}{*}{\tabincell{c}{Configuration}}    & \multicolumn{5}{c}{Considered Processes}
			& \multirow{2}{*}{\tabincell{c}{Experimental\\Validation}}  & \multirow{2}{*}{Controller}  \\
			\cline{3-7}
			& & \tabincell{c}{Load\\control}            & \tabincell{c}{Lye-gas\\ separation} & \tabincell{c}{Lye\\circulation}  & \tabincell{c}{Thermal\\control}   & \tabincell{c}{HTO\\ impurity}                          &    \\                        \hline
			\tabincell{c}{Qi et al. 2023 \cite{qi2023thermal,qi2023design}}       & 1-in-1   & $\times$      & \checkmark          & \checkmark        & \checkmark  & $\times$ 	&  \tabincell{c}{1-in-1 system} 	& MPC/PID \\
			\tabincell{c}{Qi et al. 2021 \cite{qi2021pressure}}                   & 1-in-1   & $\times$    & \checkmark          & \checkmark        & $\times$  & \checkmark  &  \tabincell{c}{1-in-1 system}  & MPC \\
			\tabincell{c}{Qiu et al. 2024 \cite{qiu2024nonlinear}}                & Multiple 1-in-1   & \checkmark    & \checkmark         & \checkmark        & \checkmark  & \checkmark &  $\times$ & MPC \\
			\tabincell{c}{Cheng et al. 2024 \cite{cheng2024self} }                             & 1-in-1   & \checkmark    & $\times$          & $\times$       & \checkmark  & $\times$ &  \tabincell{c}{1-in-1 system}  & MPPT \\
			\tabincell{c}{Li et al. 2022 \cite{li2022active}}                     & 1-in-1   & \checkmark    & \checkmark          & \checkmark        & \checkmark  & \checkmark  &  \tabincell{c}{1-in-1 system}  & PID \\
			\tabincell{c}{Li et al. 2024 \cite{li2024optimization} }              & 1-in-1   & \checkmark    & \checkmark          & \checkmark        & \checkmark  & \checkmark  &  \tabincell{c}{1-in-1 system}  & MPC \\
			\tabincell{c}{Kang et al. 2018 \cite{kang2018optimal}}                & Multiple 1-in-1  & \checkmark    & $\times$           & $\times$        & $\times$  & $\times$     & $\times$ & \tabincell{c}{Steady-state\\optimization} \\
			\tabincell{c}{Rizwan et al. 2021 \cite{rizwan2021design}}             & \tabincell{c}{$N$-in-1 (Shared \\lye circulation)}   & \checkmark    & \checkmark          & \checkmark        & \checkmark   & $\times$   & $\times$ & \tabincell{c}{Steady-state\\optimization}  \\
			\tabincell{c}{Chen et al. 2024 \cite{chen2024enhancing}}              & \tabincell{c}{$N$-in-1 (Shared \\lye circulation)}   & \checkmark    & \checkmark          & \checkmark        & \checkmark   & \checkmark   & $\times$ & PSO-MPC  \\
			\tabincell{c}{Shi et al. 2023 \cite{shi2023plant}}                   & \tabincell{c}{$N$-in-1 (Shared \\lye circulation)}   & \checkmark    & \checkmark          & \checkmark        & \checkmark   & $\times$   & $\times$ & DMC/MPC  \\	
			\tabincell{c}{Zheng et al. 2023 \cite{zheng2022optimal}}                   & \tabincell{c}{$N$-in-1 (Shared \\ lye-gas  separation \\and lye circulation)}  & \checkmark    & $\times$          & $\times$        & \checkmark   & $\times$   & $\times$ & PID/Rule-based  \\
			\tabincell{c}{Liang et al. 2024 \cite{liang2024large}}                 & \tabincell{c}{$N$-in-1 (Shared \\ lye-gas  separation \\and lye circulation)}   & \checkmark    & $\times$          & $\times$        & $\times$   & $\times$   & $\times$ & Rule-based  \\
			\tabincell{c}{Li et al. 2023 \cite{li2023exploration}}                   & \tabincell{c}{$N$-in-1 (Shared \\ lye-gas  separation \\and lye circulation)}   & \checkmark    & $\times$         & \checkmark        & \checkmark   & $\times$   & $\times$ & Rule-based  \\
			\tabincell{c}{Hu et al. 2025 \cite{hu2025operation}}                   & \tabincell{c}{$N$-in-1 (Shared \\ lye-gas  separation \\and lye circulation)}   & \checkmark    & \checkmark         & \checkmark        & \checkmark   &\checkmark  & $\times$ &
			\tabincell{c}{Rule-based}  \\
			\tabincell{c}{Guan et al. 2026 \cite{guan2026region}}                   & \tabincell{c}{$N$-in-1 (Shared \\ lye-gas  separation \\and lye circulation)}   & \checkmark    & \checkmark         & \checkmark        & \checkmark   &\checkmark  & $\times$ &
			\tabincell{c}{Two-layer\\region-based}  \\
			\hline			
			\bf This work$^{1}$     & \tabincell{c}{$N$-in-1 (Shared \\ lye-gas  separation \\and lye circulation)}   & \checkmark    & \checkmark          & \checkmark        & \checkmark   & \checkmark     &  \tabincell{c}{4,000 Nm$^3$/h-rated\\4-in-1 system} & {\color{black}MIQP-based MPC}
			\\
			\hline\hline
			\multicolumn{4}{l}{\makecell{$^{1}$ An earlier preprint version of this paper was posted as \cite{qiu2025dynamic}}}\\
		\end{tabular}
	\end{adjustbox}
\end{table}

Despite extensive research on 1-in-1 systems, studies on $N$-in-1 configurations remain limited. Existing work has mainly focused on systems sharing only the lye circulation loop, referred to as \emph{weakly coupled} systems \cite{chen2024enhancing}, while systems sharing both lye circulation and gas separation, referred to as \emph{strongly coupled} systems \cite{chen2024enhancing}, have received less attention.

For weakly coupled systems, Rizwan et al. \cite{rizwan2021design} investigated thermal management and capacity optimization and proposed a steady-state optimization method. Chen et al. \cite{chen2024enhancing} developed a model predictive control (MPC) framework integrating thermal, flow, and pressure regulation, while Shi et al. \cite{shi2023plant} proposed a dual-layer RTO-PID temperature controller. However, weakly coupled systems are less common in large-scale industrial applications. In practice, strongly coupled configurations are more widely adopted because they offer greater cost reduction \cite{meng2024advantages,zhangalkaline,zhai2024review}, although they also introduce stronger coupling effects and higher control complexity.


For strongly coupled systems, Zheng et al. \cite{zheng2022optimal} studied inter-stack load allocation under wind power conditions, although the rule-based control strategy employed cannot guarantee optimality. Other studies \cite{zhang2014wind,liang2024large,li2023exploration} focused on startup and shutdown scheduling using rule-based methods, but did not consider dynamic temperature or HTO impurity behavior.

{
\color{black}
Recently, preliminary modeling frameworks for strongly coupled $N$-in-1 AWE systems have begun to emerge. An earlier preprint by the authors \cite{qiu2025dynamic} proposed an initial state-space model for such systems. Based on this model, subsequent studies \cite{hu2025operation,guan2026region} further investigated rule-based and region-based control strategies. However, these studies mainly focused on theoretical analysis and still lack validation using industrial-scale experimental data. In addition, coordinated control and quantitative comparisons between $N$-in-1 and 1-in-1 configurations remain insufficiently studied.
}

Table \ref{tab:literature} summarizes the related literature. Compared with 1-in-1 systems, $N$-in-1 configurations introduce stronger coupling among stacks through the shared BoP system. This coupling affects temperature regulation, HTO impurity control, and operational flexibility under fluctuating renewable power inputs \cite{rizwan2021design}. Nevertheless, most existing studies on production scheduling \cite{wang2024water,varela2021modeling,qiu2023extended,zeng2024scheduling}, energy management \cite{li2024two,qiu2024flexibility,zhu2026exploring}, and ancillary services \cite{cheng2023coordinated,guan2024frequency} still assume conventional 1-in-1 systems. {\color{black}Their applicability to $N$-in-1 configurations therefore remains unclear. This leads to two important questions:
\begin{enumerate}
	\item
	How do $N$-in-1 systems compare with 1-in-1 systems under fluctuating power inputs in terms of operational flexibility and energy conversion efficiency? Does reducing capital cost compromise operational performance?
	
	\item
	How can controllers for $N$-in-1 AWE systems be designed to match the performance of 1-in-1 systems despite increased complexity?
\end{enumerate}
}


To address these questions, this study investigates $N$-in-1 AWE systems in which multiple stacks share both gas separators and the lye circulation loop. The main contributions are summarized as follows:
\begin{enumerate}
	\item
	A state-space model of $N$-in-1 AWE systems is developed, {\color{black}capturing the effects of shared lye-gas separation and lye circulation on energy conversion efficiency, lye distribution, thermal dynamics, and HTO impurity. The model is validated using experimental data from a $4,000$ Nm$^3$/h-rated 4-in-1 system.}
		
	\item
	{\color{black}A nonlinear model predictive control (NMPC) framework and a mixed-integer quadratic programming (MIQP)-based solution method} are proposed to coordinate stack and BoP operation under variable renewable power inputs. The objective is not to compare different controller structures, such as PI or linear MPC, but to evaluate the operational flexibility and dynamic performance of $N$-in-1 and 1-in-1 configurations.
	
	\item
		Comparative simulations {\color{black}based on experimentally validated models} quantify the differences between $N$-in-1 systems with different flowsheet designs and 1-in-1 systems in terms of efficiency, flexibility, and other performance metrics.
\end{enumerate}

The remainder of this paper is organized as follows. Section \ref{sec:problem} presents the system flowsheet design. Section \ref{sec:model} develops and calibrates the state-space model. Section \ref{sec:control} introduces the controller design. Section \ref{sec:case} presents simulation results. Section \ref{sec:conclusion} concludes the paper with findings and future perspectives.


\section{Flowsheets Description}
\label{sec:problem}

\subsection{System Topology}
\label{sec:tolopogy}

\begin{figure}[tb]
	\centering
	\includegraphics[scale=0.9]{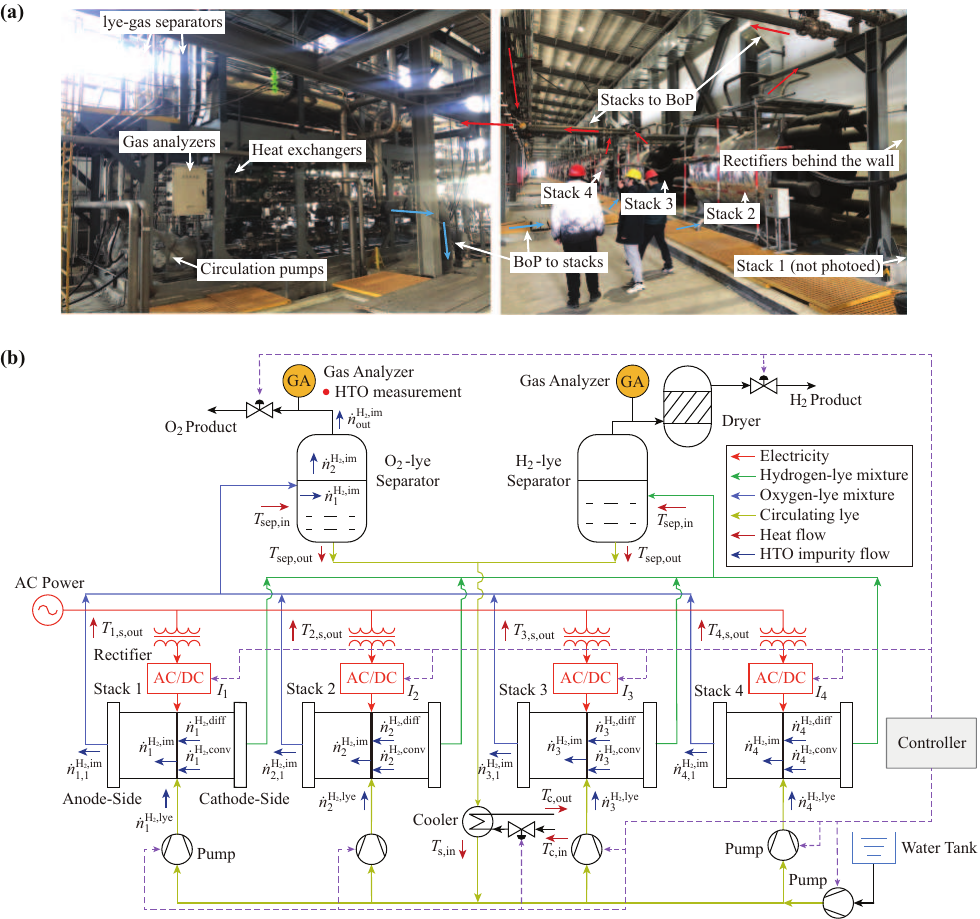}\vspace{-6pt}
	\caption{(a) {\color{black}A commercial 4,000 Nm$^3$/h-rated 4-in-1 AWE system in a large hydrogen plant}. (b) Schematic diagram of the 4-in-1 AWE system. }
	\label{fig:sys}
\end{figure}

Fig. \ref{fig:sys} shows the typical flowsheet of an $N$-in-1 AWE system. Hydrogen-lye and oxygen-lye mixtures from the stacks are mixed and sent to gas-lye separators. After separation, hydrogen and oxygen are discharged, while the lye is cooled through heat exchangers. Desalinated water is then replenished, and the lye is redistributed to the stacks. Topological variations across manufacturers are discussed in Section \ref{sec:lye}.

On the electrical side, each stack is connected to a rectifier, enabling independent control of load power. In some designs, multiple transformers are replaced by a multi-winding transformer with phase-shifted secondary windings to reduce cost and mitigate harmonic injections \cite{li2024two,zeng2024scheduling,zeng2026harmonic}.

The stacks are coupled through the shared BoP system. Their energy conversion efficiency, temperature, and HTO impurity accumulation under varying power inputs depend on coordinated control of stack power, lye flow, and cooling. To describe these interactions, a state-space model of the $N$-in-1 system is developed in Section \ref{sec:model} for controller design.

\subsection{Experimental Data Description}
\label{sec:experiment}

Experimental data, shown in Fig. \ref{fig:expdata} in Appendix A, from a 4,000 Nm$^3$/h-rated 4-in-1 AWE system in a large hydrogen plant in North China are used for model validation in Section \ref{sec:model}. The system consists of four 1,000 Nm$^3$/h-rated stacks integrated into a single BoP, with a two-pump lye circulation topology shown in Fig. \ref{fig:cir}(b).

The dataset covers 7 consecutive days, including 6 complete cycles of start-up, full-load and partial-load operation, shutdown, and warm standby.
{\color{black}
During the experiment, the rectifiers were manually controlled using stack current reference inputs, which closely overlapped with the measured stack currents shown in Fig. A1(a) due to the fast response of the power electronic devices. The circulation pumps were manually started and stopped and operated at constant speed, while the lye flow rates were regulated in open loop. The measured lye flow rates of the stacks are shown in Fig. \ref{fig:expdata}(c). The stack temperatures were regulated by controlling the inlet temperature through a PI controller, with the reference inlet temperature set to 75 $^\circ$C during production. The system pressure was maintained near 1.8 MPa during production. After shutdown, the system was depressurized through gas release, followed by a slow pressure decay caused by normal leakage before the next startup and repressurization process.}
Recorded data include system pressure, stack current and voltage, stack inlet temperatures, {\color{black}stack outlet temperatures on both the hydrogen and oxygen sides}, separator and heat exchanger temperatures, and HTO {\color{black}and OTH} impurity concentrations in the separator gas phase. The sampling interval was $1$ s.


Note that due to safety, production continuity, and cost considerations, the proposed controller was not implemented on the full-scale 4,000 Nm$^3$/h-rated system. Instead, the system dynamics were recorded experimentally to validate the model, while controller performance was evaluated through detailed simulations in Section \ref{sec:case}.

\section{State-Space Models of the $N$-in-1 AWE System}
\label{sec:model}

Compared with the 1-in-1 system, the $N$-in-1 system exhibits more complex dynamics due to coupling of heat and mass flows among stacks. Based on classical work \cite{ulleberg2003modeling} and our prior studies on 1-in-1 systems \cite{qi2023thermal,qi2023design,qi2021pressure}, a state-space model for the $N$-in-1 system is developed. The focus is on features unique to the $N$-in-1 configuration, while shared components are summarized briefly. Key variables are shown in Fig. \ref{fig:sys}(b), {\color{black} and the parameters for the 4,000 Nm$^3$/h-rated system are provided in Table \ref{tab:parameter} in Appendix B.} Without loss of generality, the following assumptions are made:

\emph{a) Pressure Dynamics:}  {\color{black} The pressure is treated as an externally specified operating condition, and is assumed to be maintained near a regulated operating value (e.g., 1.8 MPa) during production, since frequent pressure adjustment is generally avoided in practice to alleviate structural and diaphragm fatigue.}

\emph{b) Water Consumption and Replenishment:} Water consumption and replenishment are much slower than electrical and thermal dynamics and are neglected.

\emph{c) Symmetry Between Stacks:} Like most commercial $N$-in-1 systems, stacks have identical geometric parameters, and asymmetric designs are not considered.

\subsection{Electrochemical and Production Models}
\label{sec:production}

Hydrogen production is described by a semi-empirical electrochemical model \cite{sanchez2018semi}. For each stack ($i=1,\ldots,N$), the cell voltage follows
\begin{align}
  U_i^{\rm{cell}} = U_i^{\rm{rev}} + \left(r_{i,1} +r_{i,2} T_{i,\text{s}} +r_{i,3} \rho \right) I_i
  + s_i \,{\rm{log}} \left[ \left( t_{i,1} +\frac{t_{i,2}}{T_{i,\text{s}}^2} +\frac{t_{i,3}} {T_{i,\text{s}}^2} \right) I_i +1\right], \label{eq:UItotal}
\end{align}
where $U_i^{{\rm{cell}}}$, $I_i$, and $T_{\text{s},i}$ are the cell voltage, current, and stack temperature; $U_i^{\rm{rev}}=1.23$ V is the reversible voltage;
$\rho$ is the pressure; $r_{i,1}$, $r_{i,2}$, $r_{i,3}$, $s_{i}$, $t_{i,1}$, $t_{i,2}$ and $t_{i,3}$ are constant parameters.

Due to stray currents, not all currents contribute to hydrogen production. The current efficiency, known as \emph{Faraday efficiency}, is modeled by
\begin{align}
  \eta_i^{\rm{cell}} = \frac{\left(0.1 I_i\right)^2}{f_{i,1} + \left( 0.1 I_i \right)^2} f_{i,2}, \label{eq:faraday}
\end{align}
where $f_{i,1}$ and $f_{i,2}$ are coefficients related to temperature and pressure.

The hydrogen and oxygen production rates and power consumption of each stack are
\begin{align}
  \dot{n}_{i}^{\text{H}_2,\text{prod}} &= { \eta_i^{\rm{cell}} N^{\rm{cell}} I_i}/{(2F)}, \label{eq:oxygenflow} \\
  \dot{n}_{i}^{\text{O}_2,\text{prod}} &= { \eta_i^{\rm{cell}} N^{\rm{cell}} I_i}/{(4F)}, \label{eq:hydrogenflow} \\
  P_i^{\rm{ele}} &= N^{\rm{cell}} U_i^{\rm{cell}} I_i,  \label{eq:power}
\end{align}
where $F=96,485$ C/mol is the Faraday's constant; $N^\text{cell}$ is the number of cells in a stack.

Note that the energy efficiency depends on stack temperature, showing a nonlinear coupling with heat and mass transfer discussed in Section \ref{sec:thermal}. We later address them in Section \ref{sec:simplification}.
Parameter variations due to degradation can be handled by online estimation \cite{qiu2023dynamic}. The electrochemical model of individual stacks is the same as the 1-in-1 system and will not be elaborated.

\subsection{Lye Circulation, Heat Transfer and HTO Impurity Accumulation Models}
\label{sec:process}

\subsubsection{Lye Circulation and Inter-Stack Distribution}
\label{sec:lye}

\begin{figure}[t]
	\centering
	\includegraphics[scale=0.9]{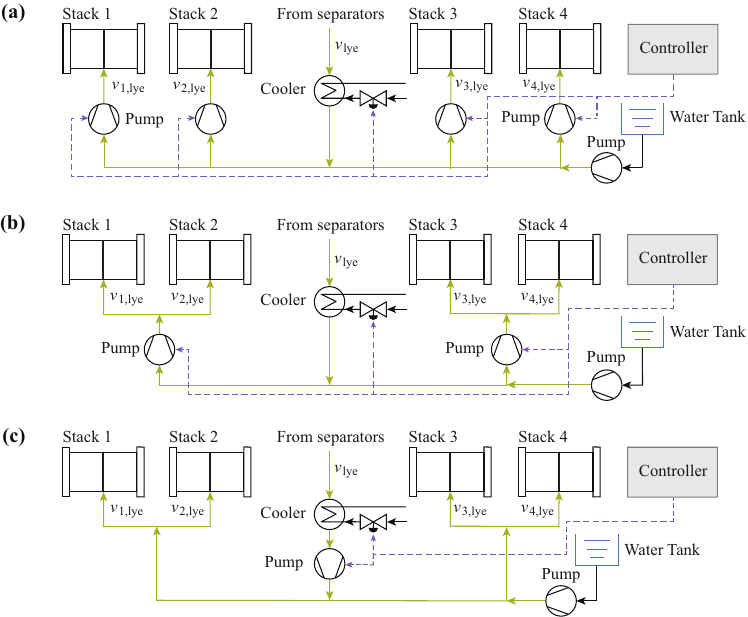}\vspace{-6pt}
	\caption{Topologies for inter-stack lye flow distribution. (a) 4 lye circulation pumps. (b) 2 lye circulation pumps. (c) 1 lye circulation pump.}
	\label{fig:cir}
\end{figure}

Fig. \ref{fig:cir} shows typical lye circulation topologies in commercial $N$-in-1 systems. Some designs use independent pumps for each stack as shown in Fig. \ref{fig:cir}(a), while others such as Figs. \ref{fig:cir}(b) and \ref{fig:cir}(c) share pumps among multiple stacks. In shared-pump systems, lye flow distribution depends on operating conditions. Using Fig. \ref{fig:cir}(c) as an example, the lye flow distribution model is derived as follows.

{\color{black}
Because lye circulation in AWE stacks is mainly designed for gas removal and temperature regulation rather than high-speed transport, the flow velocity is relatively low. Under typical operating conditions, the corresponding Reynolds number remains within the laminar or transitional regime. Therefore, the dominant pressure-flow relationship can be approximated using a Poiseuille-type hydraulic resistance model \cite{david2020dynamic}:}
\begin{align}
  \Delta p_i = v_{i,\text{mix}}^{\text{ca/an}} R_i^{\text{ca/an}}, \label{eq:poiseuille}
\end{align}
\noindent
where $\Delta p_i$ is the pressure drop across the $i$th stack; $v_{i,\text{mix}}^{\text{ca/an}}$ is the flow rate of cathode/anode-side lye-gas mixture; and $R_i^{\text{ca/an}}$ is the flow resistance.

{\color{black}
Note that under turbulent flow, the pressure drop-velocity relation (\ref{eq:poiseuille}) can be modified by including quadratic and higher-order terms. However, under the lye flow velocity in industrial AWE systems, nonlinear terms are relatively small, as experimentally verified in \cite{deiters2025experimentally}. Therefore, this work omits nonlinear terms and retains only the linear term in the pressure drop-velocity relation.
}

As stacks are connected at both ends, their pressure drops are identical. Therefore, the flow rates are inversely proportional to flow resistance:
\begin{align}
  \color{black}
   \frac{v_{i,\text{mix}}^{\text{ca/an}}}{v_{j,\text{mix}}^{\text{ca/an}}} = \frac{R_j^{\text{ca/an}}}{R_i^{\text{ca/an}}}. \label{eq:lyeallocation}
\end{align}

{\color{black}Under the above low-Reynolds-number assumption}, $R_i^{\text{ca/an}}$ can be approximated as
\begin{align}
  R_i^{\text{ca/an}} = \frac{128 \mu_{i,\text{mix}}^{\text{ca/an}} l }{\pi d^4}, \label{eq:flowresis}
\end{align}
where $l$ and $d$ are the equivalent length and diameter of the
flow paths; and $\mu_{i,\text{mix}}^{\text{ca/an}}$ is the dynamic viscosity of the lye-gas mixture, which can be approximated as a volume-weighted average of the dynamic viscosities of the lye $\mu_{\text{lye}}$ and the gas $\mu_{\text{H}_2/\text{O}_2}$ \cite{david2020dynamic}:
\begin{align}
  \mu_{i,\text{mix}}^{\text{ca/an}} = \frac{ V_{i,\text{gas}}^{\text{ca/an}} \mu_{\text{H}_2/\text{O}_2} + V_{i,\text{lye}}^{\text{ca/an}}   \mu_{\text{lye}}}{ V_{i,\text{gas}}^{\text{ca/an}} + V_{i,\text{lye}}^{\text{ca/an}}}
  = \alpha_{i}^{\text{ca/an}}  \mu_{\text{H}_2/\text{O}_2} +  \left( 1- \alpha_{i}^{\text{ca/an}}  \right) \mu_{\text{lye}},  \label{eq:viscos}
\end{align}
\noindent
where $V_{i,\text{lye}}^{\text{ca/an}}$ and $V_{i,\text{gas}}^{\text{ca/an}}$ are the volumes of the lye and gas in the $i$th stack; $\alpha_{i}^{\text{ca/an}}$ represents the volumetric gas ratio, calculated as
\begin{align}
  \alpha_{i}^{\text{ca/an}} =  \frac{V_{i,\text{gas}}^{\text{ca/an}}}{V_{i,\text{gas}}^{\text{ca/an}} + V_{i,\text{lye}}^{\text{ca/an}}} =  \frac{\dot{n}_{i}^{\text{H}_2/\text{O}_2,\text{prod}} R T_{\text{s},i} }{\rho v_{i,\text{mix}}^{\text{ca/an}}},
\end{align}
\noindent
where $R$ is the gas constant.

The flow rates of gas and lye $v_{i,\text{gas}}^{\text{ca/an}}$ and $v_{i,\text{lye}}^{\text{ca/an}}$ can then be extracted as
\begin{subnumcases}{\label{eq:vlye}}
  v_{i,\text{gas}}^{\text{ca/an}}  = \alpha_{i}^{\text{ca/an}} v_{i,\text{mix}}^{\text{ca/an}}, \label{eq:vlyegas} \\
  v_{i,\text{lye}}^{\text{ca/an}}  =  \left( 1- \alpha_{i}^{\text{ca/an}} \right) v_{i,\text{mix}}^{\text{ca/an}}, \label{eq:vlyelye}
\end{subnumcases}
\noindent
which are later used for the thermal and HTO impurity models in Sections \ref{sec:thermal} and \ref{sec:hto}.

Different electrolytic power leads to varying gas ratios, affecting flow resistances and consequently the inter-stack lye flow distribution. Compared with the 4-pump topology in Fig. \ref{fig:cir}(a), the 1-pump topology in Fig. \ref{fig:cir}(c) has less flexibility in regulating lye flow. Nevertheless, with appropriate control, the system's performance
will not be adversely affected, as detailed in Section \ref{sec:comparetopo}.

\begin{remark}
  \label{remark:1}
  {\color{black} Under the above adopted linear pressure drop-velocity relation assumption,}
  considering that the viscosity of gases (around 0.9$\times$10$^{-\text{5}}$ Pa$\cdot$s for hydrogen and 2.2$\times$10$^{-\text{5}}$ Pa$\cdot$s for oxygen \cite{nistreference}) is orders of magnitude lower than that of lye (around 2.3$\times$10$^{-\text{3}}$ Pa$\cdot$s \cite{hodges2023critical}), i.e., $\mu_{\text{H}_2/\text{O}_2} \ll \mu_{\text{lye}}$, the former can be neglected in (\ref{eq:viscos}), yielding $\mu_{i,\text{mix}}^{\text{ca/an}} \approx \big( 1-\alpha_i^{\text{ca/an}}\big) \mu_{\text{lye}} $.
  Further, by comparing (\ref{eq:vlyegas}) and (\ref{eq:vlyelye}), it follows that for different $i$ and $j$, we have $v_{i,\text{lye}}^{\text{ca/an}} \approx v_{j,\text{lye}}^{\text{ca/an}}$, meaning that even if the operating points of the stacks differ, the liquid-phase flow can be considered to distribute evenly.
  {\color{black}This is experimentally verified, as shown in Fig. \ref{fig:flowcomp}.}
  Moreover, the lye flow rates on the cathode and anode sides can be assumed equal, i.e., $v_{i,\text{lye}}^{\text{ca}} \approx v_{i,\text{lye}}^{\text{an}}$. We further denote $v_{i,\text{lye}} = v_{i,\text{lye}}^{\text{ca}} + v_{i,\text{lye}}^{\text{an}}$.
\end{remark}

\subsubsection{Heat Transfer}
\label{sec:thermal}

Temperature affects both efficiency and equipment lifetime. For the 1-in-1 AWE system, our previous studies \cite{qi2023thermal}
proposed a 3rd-order state-space model to capture the heat transfer behaviors among the stack, separators, and heat exchanger. In this work, we further extend this model to $N$-in-1 AWE systems.

{\color{black}
Note that the reaction heat generated on the cathode and anode sides is different, leading to temperature gradients near the electrodes and gas diffusion layers. Nevertheless, the lye on the two sides exchanges heat through the diaphragm, manifolds, channels, and stack structure. As the lye flows toward the stack outlets, the temperatures on the hydrogen and oxygen sides gradually converge \cite{qi2023thermal}. Experimental data shown in Fig. \ref{fig:tempcomp} further confirm the similarity of the outlet temperatures on the two sides. Therefore, since this study focuses on system-level modeling and control, the outlet and separator temperatures on the hydrogen and oxygen sides are represented using the same state variables.
}

The temperature of each stack and its internal lye flow is jointly determined by electrolytic heat, lye-gas flow, and heat dissipation to the environment, modeled as
\begin{align}
  C_{i,\text{s}} \dfrac{{\text{d}}T_{i,\text{s,out}}}{\text{d}t} = Q_{i,\text{ele}} - Q_{i,\text{s,diss}} - c_{\text{lye}} v_{i,\text{lye}} \rho_{\rm{lye}} \left( T_{i,\text{s,out}} - T_{i,\text{s,in}} \right), \label{eq:heatstack}
\end{align}
where $C_{i,\text{s}}$ denotes the total heat capacity of the stack and the lye inside it, i.e., $C_{i,\text{s}} = C_{i,\text{struc,s}} + C_{i,\text{lye,s}}$, with $C_{i,\text{lye,s}} = V_{i,\text{s,lye}} \rho_\text{lye} c_\text{lye}$; $V_{i,\text{s,lye}}$ is the volume of the lye in the stack; $\rho_\text{lye} = 1,250\ \text{kg}/\text{m}^{3}$ and $c_\text{lye} = 3,300\ \text{J}/(\text{kg}\cdot{K})$  are the density and specific heat capacity of the lye; $T_{i,\text{s,in}}$/$T_{i,\text{s,out}}$ denotes the inlet/outlet lye temperature of the $i$th stack;
$Q_{i,\text{ele}}$ denotes heating flow; {\color{black}$Q_{i,\text{s,diss}}$} denotes heat dissipation to the environment.

Note that, as indicated in Fig. \ref{fig:sys}(b), the inlet temperature of the stacks is the same, while the outlet temperature varies due to differences in the temperature and load among stacks.

The heating flow $Q_{i,\text{ele}}$ comprises the heat released by the electrolytic reaction and the ohmic heat from the stray current, as
\begin{align}
  Q_{i,\text{ele}} = \eta_i^{\text{cell}} N^{\text{cell}} I_i \left( U_i^{\text{cell}} - U^{\text{th}} \right) + \left( 1 - \eta_i^{\text{cell}} \right) N^{\text{cell}} I_i U_i^{\text{cell}}, \label{eq:heatele}
\end{align}
where $U^\text{th}=1.48$ V is the thermoneutral voltage.

{\color{black}
Since the electrolyzers are installed indoors without forced ventilation directly acting on the stack surfaces, as shown in Fig. 1(a), the air velocity around the equipment is relatively low. Therefore, the heat dissipation $Q_{i,\text{s,diss}}$ from the $i$th stack to the environment is approximated as the sum of natural convection $Q_{i,\text{s,conv}}$ and thermal radiation $Q_{i,\text{s,rad}}$, which follow}
\begin{align}
  Q_{i,\text{s,diss}} = Q_{i,\text{s,conv}} + {\color{black}Q_{i,\text{s,rad}}} = h_\text{s} A_{\text{s,diss}} ( T_{i,\text{s,out}} - T_\text{am}) + \sigma A_{\text{s,diss}} \epsilon_\text{s} ( T^4_{i,\text{s,out}} - T^4_\text{am} ), \label{eq:heatdiss}
\end{align}
where $A_{\text{s,diss}}$ is the  heat dissipation area;
$T_\text{am}$ is the ambient temperature; $\sigma$ is the Stefan-Boltzmann constant; $\epsilon_\text{s}$ denotes the emissivity of the stack; and $h_\text{s}$ is the natural convection coefficient,
following \cite{dieguez2008thermal}:
\begin{align}
  h_\text{s} = 2.51 \times 0.52 \left[ \left( {T_{i,\text{s,out}} - T_{\text{am}}}  \right) / {\phi_{s}} \right]^{0.25}, \label{eq:disscoef}
\end{align}
where $\phi_{s}$ represents the diameter of the stack.

As indicated in Fig. \ref{fig:sys}, the lye exiting the $N$ stacks is mixed before entering the separators. Assuming no heat loss in pipes, the separator inlet temperature $T_{\text{sep,in}}$ is the weighted average of the outlet lye temperatures of the $N$ stacks: 
\begin{align}
  T_{\text{sep,in}} = \dfrac{ \sum_{i=1}^N v_{i,\text{lye}} T_{i,\text{s,out}} } {v_{\text{tot,lye}}},
\end{align}
where $v_{\text{tot,lye}} = \sum_{i=1}^N v_{i,\text{lye}}$ denotes total lye flow.

For the separators, assuming the heat carried by gaseous hydrogen/oxygen is negligible, the temperature is determined by
\begin{align}
  C_{\text{sep}} \dfrac{{\text{d}}T_{\text{sep,out}}}{\text{d}t} = \frac{1}{2} c_{\text{lye}} v_{\text{tot,lye}} \rho_{\rm{lye}} \left( T_{\text{sep,in}} - T_{\text{sep,out}} \right) - Q_{\text{sep,conv}} - Q_{\text{sep,rad}}  , \label{heatseparator}
\end{align}
where $C_{\text{sep}}$ is the heat capacity of the separator and lye inside it; {\color{black}$T_{\text{sep,out}}$} is the separator outlet temperature; $Q_{\text{sep,conv}}$ and $Q_{\text{sep,rad}}$ are the convective and radiative dissipation, calculated similarly to (\ref{eq:heatdiss})--(\ref{eq:disscoef}).

The lye leaving the cathode- and anode-side separators is mixed and cooled before being circulated to the stack.  Its temperature behavior is modeled as
\begin{align}
  C_{\text{he}} \dfrac{{\text{d}}T_{\text{s,in}}}{\text{d}t} = c_{\text{lye}} v_{\text{tot,lye}} \rho_{\rm{lye}} \left( T_{\text{sep,out}} - T_{\text{s,in}} \right) - k A_\text{c} \Delta T  , \label{eq:heatexchange}
\end{align}
where $C_{\text{he}}$ is the total heat capacity of the heat exchanger and the lye it contains; 
$A_\text{c}$ represents heat exchange area; $k$ is the heat transfer coefficient;  and $\Delta T$ is the logarithmic mean temperature difference.

{\color{black}
Assuming the heat exchanger is a commonly used counterflow type, as illustrated in Fig. \ref{fig:sys}(b), the lye stream from the separator outlet to the stack inlets is treated as the hot fluid, with the heat exchanger inlet and outlet temperatures denoted by $T_{\text{sep,out}}$ and $T_{\text{s,in}}$, respectively. The cooling water is treated as the cold fluid, with inlet and outlet temperatures denoted by $T_{\text{c,in}}$ and $T_{\text{c,out}}$, respectively. Thus, $\Delta T$ is given by}
\begin{align}
	\Delta T =  \frac{ {\color{black} \left( T_{\text{sep,out}} - T_{\text{c,out}} \right) - \left( T_{\text{s,in}} - T_{\text{c,in}} \right) } }{\log  \left[ {\color{black} \left( T_{\text{sep,out}} - T_{\text{c,out}} \right) /  \left( T_{\text{s,in}} - T_{\text{c,in}} \right) } \right]},
\end{align}

{\color{black}
The inlet cooling water temperature $T_{\text{c,in}}$ is regulated by a chiller and is assumed to be constant. The outlet temperature of the cooling water, $T_{\text{c,out}}$, is then determined as follows:}
\begin{align}
  C_{\text{c}} \dfrac{{\text{d}}T_{\text{c,out}}}{\text{d}t} = c_{\text{c}} v_{\text{c}} \rho_{\text{c}} \left( T_{\text{c,in}} - T_{\text{c,out}} \right) + k A_\text{c} \Delta T  , \label{eq:heatcooling}
\end{align}
where $C_{\text{c}}$ denotes the total heat capacity of the cooling coil and the cooling water it contains;
$c_{\text{c}}=4,100$ J/(kg$\cdot$K) and $\rho_{\text{c}} = 1,000$ kg/m$^3$ are the specific heat capacity and density of the cooling water, respectively; and $v_{\text{c}}$ represents the cooling water flow rate, which is regulated by a valve.

For an $N$-in-1 system, the heat transfer model (\ref{eq:heatstack})--(\ref{eq:heatcooling}) has an order of $N+3$.
This model links temperature dynamics with control of electrolytic load, lye flow, and cooling water.
The temperatures also appear in the electrochemical model in Section \ref{sec:production} and the HTO impurity accumulation model in Section \ref{sec:hto},
showing coupling relations, which are addressed in the controller design in Section \ref{sec:control}.

\subsubsection{HTO Impurity Accumulation}
\label{sec:hto}

HTO impurity arises from gas crossover and lye mixing, affects safe operation at low loads, and must be accounted for in controller design.
Here, we extend the 3rd-order model for 1-in-1 systems proposed in our previous work \cite{qi2021pressure} to the $N$-in-1 system.

In each stack, {\color{black}the overall HTO impurity crossover from the cathode-side half-cell to the anode-side half-cell, denoted as $\dot{n}^{\text{H}_2,\text{im}}_i$, is composed of three parts, each related to lye circulation \cite{david2020dynamic}, diffusion \cite{trinke2018hydrogen}, and convection \cite{trinke2018hydrogen}}:
\begin{align}
  \dot{n}^{\text{H}_2,\text{im}}_i = \dot{n}^{\text{H}_2,\text{lye}}_i +  \dot{n}^{\text{H}_2,\text{diff}}_i +  \dot{n}^{\text{H}_2,\text{conv}}_i  , \label{eq:nim}
\end{align}
\noindent
where  {\color{black}$\dot{n}^{\rm{H_2}}_{\rm{lye}}$, $\dot{n}^{\rm{H_2}}_{\rm{diff}}$, and $\dot{n}^{\rm{H_2}}_{\rm{conv}}$} are the molar flows of HTO impurity brought in by lye circulation, diffusion, and convection, respectively.

As shown in Fig. \ref{fig:sys}, {\color{black}the circulating lye contains dissolved and unseparated hydrogen because the return streams from the hydrogen-side and oxygen-side separators are mixed. The hydrogen impurity mainly originates from the hydrogen-side separator. After mixing, the impurity concentration is approximately halved. The mixed lye is then approximately evenly distributed between the cathode-side and anode-side half-cells.} The resulting corresponding molar flow of HTO impurity into the anode-side half-cell, $\dot{n}^{\text{H}_2,\text{lye}}_i$, depends on the distribution of lye flow across the $N$ stacks (see Section \ref{sec:lye}), following \cite{david2020dynamic}
\begin{align}
  \dot{n}^{\text{H}_2,\text{lye}}_i = { {\color{black}S^{\text{H}_2}_{\text{eff}}} \rho v_{i,\text{lye}}}/{4},\label{eq:nlye}
\end{align}
{\color{black}where $S^{\text{H}_2}_{\text{eff}}$ is the effective hydrogen carrying coefficient in the circulating lye. Due to unseparated microbubbles and gas holdup, $S^{\text{H}_2}_{\text{eff}}$ can be higher than the theoretical thermodynamic equilibrium solubility and needs to be calibrated based on operation data.}

The crossover molar flows of HTO impurities due to diffusion and convection,  $\dot{n}^{\text{H}_2,\text{diff}}_i$ and $\dot{n}^{\text{H}_2,\text{conv}}_i$, are determined by Fick's law and Darcy's law \cite{trinke2018hydrogen}, respectively, as
\begin{subnumcases}{\label{eq:hto}}
  \dot{n}^{\text{H}_2,\text{diff}}_i  = A^{\text{cell}} N^{\text{cell}} \frac{D^{\text{H}_2}_{\text{eff}} \Delta c^{\rm{H_2}}}{\delta} \approx A^{\text{cell}} N^{\text{cell}} \frac{D^{\text{H}_2}_{\text{eff}} S^{\text{H}_2}_{\text{eff}} \rho}{\delta},\label{eq:ndiff} \\
  \dot{n}^{\text{H}_2,\text{conv}}_i  = A^{\text{cell}} N^{\text{cell}} \frac{K^{\text{H}_2}_{\text{eff}}}{\mu_{i,\text{lye}}} S^{\text{H}_2}_{\text{eff}} \rho \frac{\Delta \rho}{\delta},\label{eq:nconv}
\end{subnumcases}
where $D^{\rm{H_2}}_{\rm{eff}}$ is the diffusion coefficient;
$\Delta c^{\rm{H_2}}$ denotes the differential concentration of hydrogen between the cathode and anode sides; $\delta$ is the thickness of the diaphragm; $K^{\text{H}_2}_{\text{eff}}$ is the permeability of the hydrogen through the diaphragm; $\Delta \rho$ denotes the pressure difference between two sides, usually caused by fluctuations in pressure and liquid levels.

The HTO impurity is transported and accumulates in three stages, including the anode half-cell, and the liquid and gas phases  of the separator, as shown in Fig. \ref{fig:sys}.
Based on molar volume conservation, we establish a state-space model for HTO accumulation in the $N$-in-1 system, as
\begin{subnumcases}{\label{eq:Nbalance}}
  \dfrac{ \text{d} n^{\text{H}_2,\text{an}}_{i}} {\text{d}t} = \dot{n}^{\text{H}_2,\text{im}}_{i} - \dot{n}^{\text{H}_2,\text{im}}_{i,1}, \ i=1, \ldots, N, \\
  \dfrac{\text{d} n^{\text{H}_2,\text{sep}}_{\text{liq}}}{\text{d}t} = \sum\nolimits_{i=1}^{N} \dot{n}^{\text{H}_2,\text{im}}_{i,1} - \dot{n}^{\text{H}_2,\text{im}}_2,  \\
  \dfrac{\text{d} n^{\text{H}_2,\text{sep}}_{\text{gas}}}{\text{d}t} = \dot{n}^{\text{H}_2,\text{im}}_{2} - \dot{n}^{\text{H}_2,\text{im}}_{\text{out}},
\end{subnumcases}
where $n^{\text{H}_2,\text{an}}_{i}$, $ n^{\text{H}_2,\text{sep}}_{\text{liq}}$, and {\color{black}$ n^{\text{H}_2,\text{sep}}_{\text{gas}}$} represent the molar quantities of HTO impurity in the anode half-cell of the $i$th stack, and the liquid and gas phases of the separator, respectively; $ \dot{n}^{\text{H}_2,\text{im}}_{i,1}$, $\dot{n}^{\text{H}_2,\text{im}}_2$, and $\dot{n}^{\text{H}_2,\text{im}}_{\text{out}}$ are the molar flows of HTO impurity, 
which satisfies
\begin{subnumcases}{\label{eq:nbalance}}
  \dot{n}^{\text{H}_2,{\color{black}\text{im}}}_{i,1} = \frac{n^{\text{H}_2,\text{an}}_{\color{black}i} v_{i,\text{lye}} }{2 V^{\text{an}}_{{\color{black}i},\text{lye}}}, \ i=1,\ldots,N,\\
  \dot{n}^{\text{H}_2,{\color{black}\text{im}}}_2 =\frac{n^{\rm{H_2,sep}}_{\rm{liq}}}{\tau_{\rm{sep}}}, \\
  \dot{n}^{\text{H}_2,{\color{black}\text{im}}}_{\rm{out}}=  \frac{R T_{\text{sep,out}} n^{\text{H}_2,\text{sep}}_{\text{gas}} \sum\nolimits_{i=1}^{N} \dot{n}^{\text{O}_2,\text{prod}}_i} {\rho V^{{\color{black}\text{sep}}}_{{\color{black}\text{gas}}}}, \label{eq:nbalance3}
\end{subnumcases}
where the anode-side half-cell lye volume {\color{black}$V_{i,\rm{lye}}^{\rm{an}}$} is determined by (\ref{eq:vlye}); $\tau_{\rm{sep}}$ is the separation time constant; $V^{{\color{black}\text{sep}}}_{{\color{black}\text{gas}}}$ is the gas-phase volume of the separator.

The HTO impurity concentration in the gas phase of the oxygen-side separator, which is measured and taken as the safety metric \cite{david2020dynamic,trinke2018hydrogen,qi2021pressure}, follows:
\begin{align}
  \text{HTO} = \frac{n^{\text{H}_2,\text{sep}}_{\text{gas}} R T_{\text{sep,out}}} {\rho V^{{\color{black}\text{sep}}}_{{\color{black}\text{gas}}}}. \label{eq:defhto}
\end{align}

Compared to the 1-in-1 system model in \cite{qi2021pressure}, the model order of the $N$-in-1 system increases to $N$+2. The accumulation process depends not only on the overall power but also on inter-stack load allocation and lye flow control, significantly increasing the complexity.
The HTO impurity accumulation process is also integrated into the controller to accommodate the varying energy input;
see Section \ref{sec:control} for details.

{\color{black}
The oxygen-to-hydrogen (OTH) impurity concentration in the hydrogen-side separator is much lower than the HTO concentration \cite{emam2024review,qi2021pressure,david2020dynamic}. As shown in the experimental data in Fig. \ref{fig:expdata}(f), the OTH concentration is mostly below 0.1\%. While HTO impurity is a key safety constraint, OTH impurity is less critical. Therefore, we leave OTH impurity dynamics for future refined studies.	
}

{\color{black}
Note that although this work adopts a quasi-stationary pressure assumption, rapid pressurization or depressurization processes may influence transient HTO dynamics. In particular, gas holdup, venting behavior, separator liquid-level variations, and pressure-differential dynamics may affect the transient evolution of impurity concentration during fast pressure transients. Although these effects mainly influence transient impurity responses, their influence on the long-timescale operating equilibrium is comparatively smaller under pressure-regulated operating conditions. Incorporating coupled pressure dynamics into the HTO model would further improve the fidelity of transient safety analysis, which merits future study.
}

\begin{table}[tb]\scriptsize\centering
	\renewcommand{\arraystretch}{1.35}
	\caption{Summary of the proposed state-space model of the $N$-in-1 AWE system}\vspace{6pt}
	\label{tab:model}
	\begin{tabular}{ ccccc }\hline\hline
		Submodel              & \tabincell{c}{Electrochemical\\ and production}   & \tabincell{c}{Lye flow\\distribution} & \tabincell{c}{Heat transfer}                & \tabincell{c}{HTO impurity accumulation} \\ \hline
		State variables       & /                                                 & /                                                 & $T_{\text{s,in}}$, $\big\{ T_{i,\text{s,out}} \big\}_{i=1}^N$, $T_{\text{sep,out}}$, $T_{\text{c,out}}$                                              & $\big\{ n_i^{\text{H}_2,\text{an}} \big\}_{i=1}^N$,  $ n_{\text{liq}}^{\text{H}_2,\text{sep}}$,  $ n_{\text{gas}}^{\text{H}_2,\text{sep}} $                                         \\
		Order of model        & $0$                                               & $0$                                               & $N+3$                                         & $N+2$ \\
		Formulation           & (\ref{eq:UItotal})--(\ref{eq:power})              & (\ref{eq:poiseuille})--(\ref{eq:vlye})            & (\ref{eq:heatstack})--(\ref{eq:heatcooling})  & (\ref{eq:nim})--(\ref{eq:defhto})\\
		Control variables     &  \multicolumn{4}{c}{ Electrolytic currents $\big\{ I_i \big\}_{i=1}^N$, lye flow rates $\big\{ v_{i,\text{lye}} \big\}_{i=1}^N$, and cooling water flow rates $v_\text{c}$}        \\
		\hline\hline
	\end{tabular}
\end{table}

\begin{figure}[H]
	\centering
	\includegraphics[scale=0.90]{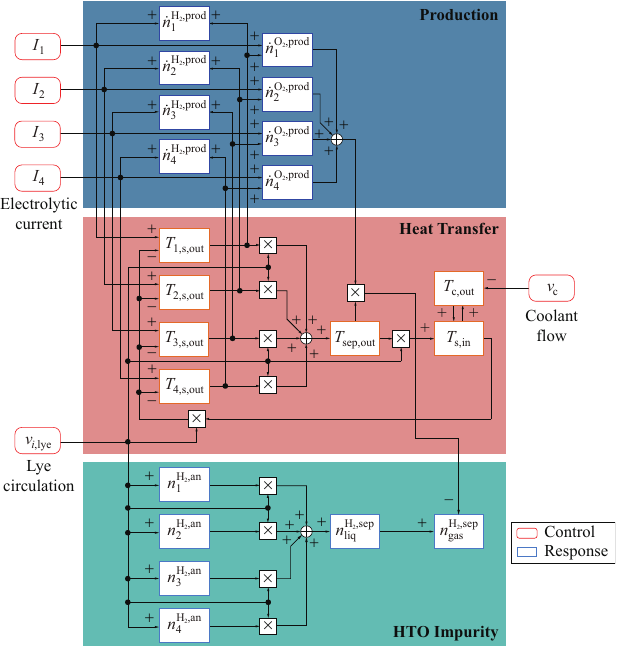}\vspace{-6pt}
	\caption{Interrelationships between state variables, control, and responses in the $N$-in-1 AWE system with the 1-pump topology in Fig. \ref{fig:cir}(c).}
	\label{fig:state}
\end{figure}

\subsection{Model Summary}
\label{sec:modelsummary}

Summarizing the submodels established in Sections \ref{sec:production} and \ref{sec:process}, Table \ref{tab:model} presents an overview of the state-space model for the $N$-in-1 AWE system. Note that the electrochemical and inter-stack lye distribution submodels are formulated by algebraic equations. The complete state-space model has an order of $2N+5$.

For easy understanding, Fig. \ref{fig:state} illustrates the coupling among hydrogen production, temperature, and HTO impurity; see simulations in Section \ref{sec:openloop}.
Compared to the 1-in-1 system, the $N$-in-1 system introduces additional control freedom through current and lye flow allocation, making optimal control more challenging.

The state and control variables can be summarized in vector form:
\begin{gather}
  \bm{x} \triangleq \big[ T_{\text{s,in}}, T_{1,\text{s,out}}, \ldots, T_{N,\text{s,out}}, T_{\text{sep,out}}, T_{\text{c,out}}, n_1^{\text{H}_2,\text{an}}, \ldots, n_N^{\text{H}_2,\text{an}}, n_{\text{liq}}^{\text{H}_2,\text{sep}}  , n_{\text{gas}}^{\text{H}_2,\text{sep}} \big]^{\text{T}}, \label{eq:state} \\
  \bm{u}  \triangleq \big[  I_1, \ldots, I_N, v_{1,\text{lye}}, \ldots, v_{N,\text{lye}}, v_\text{c} \big]^{\text{T}}, \label{eq:control}
\end{gather}
\noindent
and the state equations can therefore be compactly denoted as
\begin{align}
  {\text{d}\bm{x}}/{\text{d}t} = \bm{h}(\bm{x}, \bm{u}). \label{eq:equstate}
\end{align}

\subsection{Model Validation Using Experimental Data}
\label{sec:validation}

\begin{figure}[tb]
	\centering
	\includegraphics[scale=0.94]{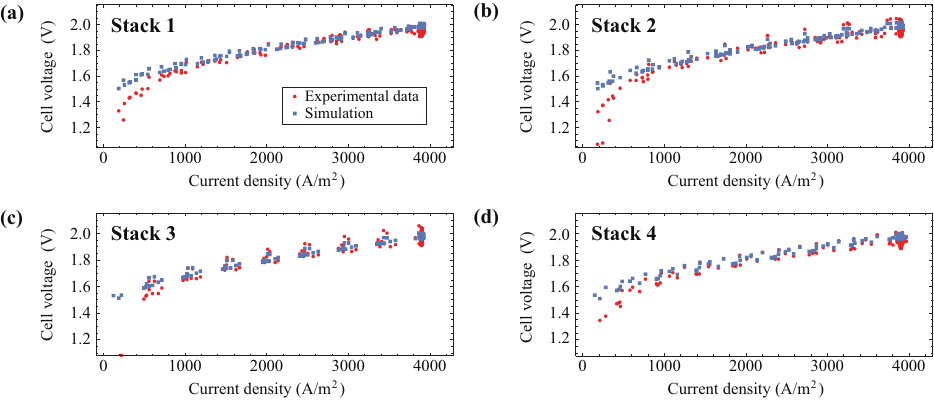}\vspace{-6pt}
	\caption{{\color{black}Validation of the cell voltage model against the experimental data of the 4,000 Nm$^3$/h-rated 4-in-1 electrolysis system. (a) Stack 1. (b) Stack 2. (c) Stack 3. (d) Stack 4. }}
	\label{fig:voltvalidation}
\end{figure}

\begin{figure}[tb]
	\centering
	\includegraphics[scale=0.94]{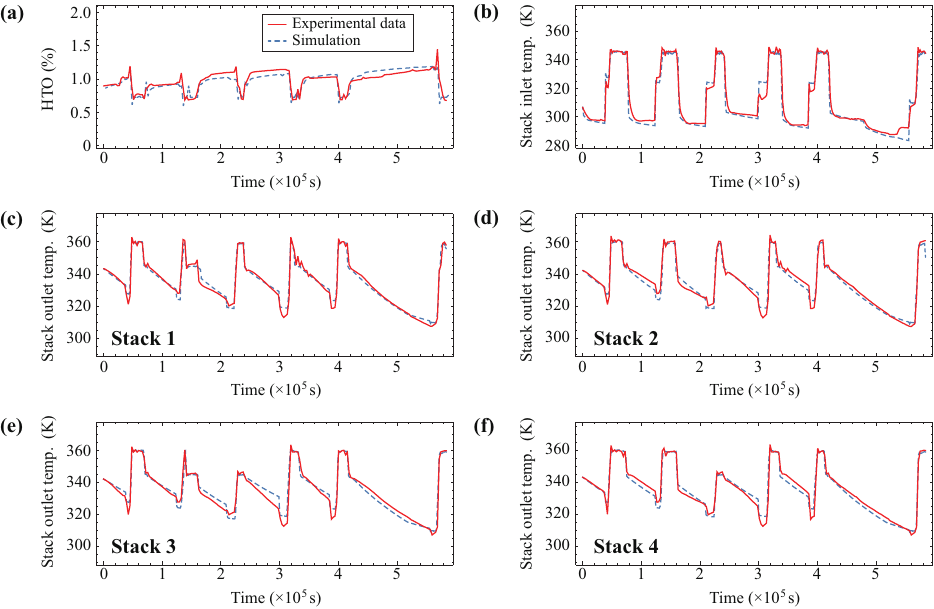}\vspace{-6pt}
	\caption{{\color{black}Validation of the proposed state-space model against the experimental data of the 4,000 Nm$^3$/h-rated 4-in-1 electrolysis system. (a) HTO impurity in the gas phase of the oxygen-lye separator. (b) Inlet temperature of the stacks. (c)--(f) Outlet temperatures of the stacks. }}
	\label{fig:validation}
\end{figure}

The experimental data described in Section \ref{sec:experiment} and Appendix A are used to validate the model and calibrate its parameters. {\color{black}First, after fixing available parameters such as geometric dimensions and rectifier ratings from factory datasheets, electrochemical, heat transfer, and mass transfer model parameters are identified using the estimation method proposed in \cite{qiu2023dynamic}, with results listed in Table \ref{tab:parameter} in Appendix B.}
The calibrated model is then simulated.  {\color{black}The measured trajectories of control inputs and operating conditions (currents, pressure, lye flow rates, etc., shown in Fig. \ref{fig:expdata}) are directly used as simulation inputs.}
{\color{black} The validation of cell voltage output is shown in Fig. \ref{fig:voltvalidation}, and the results of temperature and HTO impurity dynamics are compared with measurements in Fig. \ref{fig:validation}.}

{\color{black}
As shown in Fig. \ref{fig:voltvalidation}, the semi-empirical electrochemical model agrees well with the experimental data when the current density exceeds 1,000 A/m$^2$. The measured data fluctuate around a smooth trend because the operating temperature and pressure vary during operation. Since the model is intended to describe the electrochemical behavior after stack activation, larger deviations appear in the low-load region. However, these deviations are acceptable because the voltage model has limited influence on the dynamic heat and mass transfer processes under low-load conditions. In practice, low-load operation is also generally avoided due to low Faraday efficiency and increased HTO impurity concentration. Overall, the root-mean-square errors (RMSEs) of the cell voltage for stacks 1--4 are 0.02819, 0.03046, 0.03294, and 0.02997 V, respectively, indicating sufficient accuracy for system-level modeling and control.
}

As shown in Fig. \ref{fig:validation}, the heat and mass transfer models also closely match the field data over the full 7-day operation period (approximately $6\times10^{5}$ s). The RMSEs of stack inlet and outlet temperatures are $4.001$ K and $3.682$ K, respectively, while the RMSE of HTO impurity concentration is 0.1387\%. These errors are small given the environmental uncertainty in long-duration field experiments.

{\color{black}
To further validate the generalization ability of the model, we also perform parameter calibration and validation using split datasets. The 7-day dataset is divided chronologically into two parts. The data from the first three startup/shutdown cycles are used for parameter calibration, while the remaining data are used for independent validation under unseen operating conditions. The calibrated model is then directly applied to the validation dataset. The resulting RMSEs for the calibration dataset are 3.517 K, 2.151 K, and 0.132\% for stack inlet temperature, stack outlet temperature, and HTO concentration, respectively; while for the validation dataset the corresponding RMSEs are 3.820 K, 4.020 K, and 0.148\%. The results indicate that the proposed model maintains good predictive capability across different operating periods.
}

{\color{black}Considering that the 7-day dataset covers complete operating cycles, including start-up, full- and partial-load operation, standby, and shutdown, the results indicate that the model captures system dynamics with sufficient accuracy, and we use the 7-day-dataset calibrated parameters to support controller design in Section \ref{sec:control} and performance evaluation in Section \ref{sec:case}.}

\section{Controller Design}
\label{sec:control}

\subsection{Overview of Control}

Section \ref{sec:modelsummary} shows that the $N$-in-1 AWE system is a nonlinear multi-input multi-output (MIMO) system, where safety and efficiency must be balanced under fluctuating power inputs.
Conventional PID controllers \cite{qi2023design,qi2021pressure} and linear MPC based on local linearization are inadequate for such dynamics \cite{baldea2014integrated}.
To address this, we develop an NMPC-based framework.

As in standard MPC, a control horizon $N^\text{h}$ and time step $\Delta t$ are defined.
Because temperature and HTO accumulation evolve over tens of minutes, a relatively long horizon is required. To limit computational burden while preserving accuracy, we discretize (\ref{eq:equstate}) using the trapezoidal scheme to ensure accuracy and numerical stability:
\begin{align}
  \frac{\bm{x}(k+1) - \bm{x}(k) }{\Delta t} = \frac{\bm{h}\big( \bm{x}(k+1), \bm{u}(k+1) \big) + \bm{h}\big( \bm{x}(k), \bm{u}(k) \big)}{2}, \ k = 0, \ldots, N^{\text{h}}-1, \label{eq:diststate}
\end{align}

The control objectives and constraints are established as follows.

\subsection{Control Objective}
\label{sec:obj}

The controller aims to maximize hydrogen production while tracking the available power. It also regulates temperature to limit thermal stress and avoids frequent changes in lye and cooling flows. {\color{black}The objective function is defined in terms of minimizing the cost, as follows:}
\begin{align}
  f\big( \bm{x}(1), \ldots, \bm{x}(N^\text{h}), & \bm{u}(0), \bm{u}(1), \ldots, \bm{u}(N^\text{h}) \big) = - \lambda^{\text{prod}} \Delta t \sum\nolimits_{k=0}^{N^\text{h}} \sum\nolimits_{i=1}^{N} \dot{n}_{i}^{\text{H}_2,\text{prod}}  \nonumber \\
  + & \lambda^{\text{track}} \sum\nolimits_{k=0}^{N^\text{h}} \Big( P_{\text{tot}}^{\text{ref}}(k) - \sum\nolimits_{i=1}^{N} P_{i}^{\text{ele}}(k) \Big)^2  \nonumber \\
  + & \lambda^{\text{temp}} \sum\nolimits_{k=0}^{N^\text{h}} \sum\nolimits_{i=1}^{N} \Big( T_{i,\text{s,out}}(k) -  T_{\text{s,out}}^{\text{ref}} \Big)^2  \nonumber \\
  + & \lambda^{\text{I}} \sum\nolimits_{k=0}^{N^\text{h}-1} \sum\nolimits_{i=1}^{N} \Big( I_{i}(k+1) -  I_{i}(k) \Big)^2  \nonumber \\
  + & \lambda^{\text{lye}} \sum\nolimits_{k=0}^{N^\text{h}} \sum\nolimits_{i=1}^{N} \Big( v_{i,\text{lye}}(k) -  v_{i,\text{lye}}^{\text{0}} \Big)^2  \nonumber \\
  + & \lambda^{\text{c}} \sum\nolimits_{k=0}^{N^\text{h}}  \Big( v_{\text{c}}(k) -  v_{\text{c}}^{0} \Big)^2,
   \label{eq:obj}
\end{align}
\noindent
where $P_{\text{tot}}^{\text{ref}}$ is the reference for total power consumption; $T_{\text{s,out}}^{\text{ref}}$ is the temperature reference for the $i$th stack; $\lambda^{\text{prod}}$, $\lambda^{\text{track}}$, $\lambda^{\text{temp}}$, $\lambda^{\text{I}}$, $\lambda^{\text{lye}}$, and $\lambda^{\text{c}}$  are weight factors.
In simulations, we find that these weights can be changed with a relatively large range without significantly impacting the control performance. Hence, to save space, we will not discuss the selection of these weights in detail.

For power tracking, i.e., the second row in (\ref{eq:obj}), three application scenarios are considered:
\begin{enumerate}
  \item{Tracking renewable energy.} Forecasts of wind and solar power are used as the power reference $P_{\text{tot}}^{\text{ref}}(k)$, and the forecast error is alleviated through the rolling implementation of NMPC.
      Techniques such as robust control \cite{li2024two}
      can be used to improve performance under uncertainty.

  \item{Participating in peak-shaving for the power grid: Power references are issued by grid operators in advance, which is a deterministic time series.}

  \item{ Tracking automatic generation control (AGC) signals}: Due to strong uncertainty and temporal correlation, this requires stochastic MPC (SMPC)  \cite{qiu2021continuous,qiu2020stochastic}. Since this paper focuses on the AWE system, deterministic power references are considered, leaving SMPC for future research.
\end{enumerate}

\subsection{Constraints}
\label{eq:constraints}

At each time step $(k = 1,\ldots, N^\text{h})$ and for each stack $(i=1,\ldots,N)$, the following constraints must be satisfied to ensure operational feasibility and safety.

\subsubsection{DC Power and Cell Voltage Constraints}

The rectifier that provides DC power has current and power limits,
expressed as:
\begin{gather}
  0 \le I_i(k) \le \overline{I}, \ \text{and} \  
  0 \le P^\text{ele}_i(k) \le \overline{P}^\text{ele}, \label{eq:limitdcpower}
\end{gather}
\noindent
where $\overline{I}$ and $\overline{P}^\text{ele}_i$ are the current and power limits that the rectifier can provide, which is usually $1.2$ to $1.4$ times the rated power of the stack.
Moreover, to avoid excessive stress on electrodes and catalysts, the cell voltage is limited:
\begin{align}
  U_i^\text{cell}(k) \le \overline{U}^\text{cell}, \label{eq:limitcellvolt}
\end{align}
where $\overline{U}^\text{cell}$ is typically set at $2.1$ to $2.3$ V.

Note that some studies \cite{rizwan2021design,chen2024enhancing,shi2023plant,li2023exploration,liang2024large} and \cite{zheng2022optimal} set a non-zero lower power limit to avoid HTO over-limit under low loads. However, HTO accumulation is a dynamic process. Setting non-zero lower current or power limits is prone to being conservative by excluding short-term low-load operations \cite{qiu2023extended}.
We integrate HTO dynamics into the controller to exploit flexibility, eliminating the need for non-zero lower power limits.

\subsubsection{Hydrogen Production and Ramp Constraints}

In order to avoid sudden changes in pressure, the gas-liquid ratio in the stacks, separator liquid levels, etc., which could induce structural stress, the hydrogen production rate is subject to ramp constraints:
\begin{align}
  r^{\text{H}_2,\text{prod,down}} \le \frac{\dot{n}_{i}^{\text{H}_2,\text{prod}}(k) - \dot{n}_{i}^{\text{H}_2,\text{prod}}(k-1)}{\Delta t} \le r^{\text{H}_2,\text{prod,up}} , \label{eq:limitramp}
\end{align}
\noindent
where $r^{\text{H}_2,\text{prod,up}}$ and $r^{\text{H}_2,\text{prod,down}}$  are the upward and downward ramp limits, respectively.

\subsubsection{Temperature and HTO Constraints}

To ensure the temperature and HTO impurity concentration stay within safe limits, the following constraints are imposed:
\begin{gather}
  T_{i,\text{s,out}}(k) \le \overline{T}, \text{and}\ \text{HTO}(k) \le 2\%.  \label{eq:limittemp}
\end{gather}

\subsubsection{Lye and Coolant Flow Constraints}

Due to the limitations of pump load, valve opening range, cold water inlet pressure, etc., the lye flow and cooling water flow must also be controlled within the feasible range, as:
\begin{gather}
   \underline{v}_{\text{lye}} \le v_{i,\text{lye}}(k) \le \overline{v}_{\text{lye}}, \ \text{and} \  
   \underline{v}_{\text{c}} \le v_{\text{c}}(k) \le \overline{v}_{\text{c}},  \label{eq:limitcool}
\end{gather}
\noindent
where $\underline{v}_{\text{lye}}$, $\overline{v}_{\text{lye}}$ and $\underline{v}_{\text{c}}$, $\overline{v}_{\text{c}}$ denote the respective lower and upper limits of lye and coolant flow rates.

\subsection{Simplification of Nonlinear Components}
\label{sec:simplification}

The complete state-space model developed in Section \ref{sec:model} is highly nonlinear and non-convex, making direct integration into NMPC computationally infeasible. To address this, the following simplifications are made.

\subsubsection{Polyhedral Approximation of Production Function and Electrolytic Heat}
\label{sec:polynomial}

From the models (\ref{eq:UItotal})--(\ref{eq:power}) in Section \ref{sec:production}, it can be seen that the current, temperature, electrical power, and hydrogen flow show strong nonlinearity. Fortunately, if the relationship between {\color{black} stack} power and gas production is relaxed into the form of inequality constraints:
\begin{gather}
	\dot{n}_{i}^{\text{H}_2,\text{prod}} (I_i,T_{i,\text{s}}) \le {\color{black} \dfrac{ \eta_{i}^{\text{cell}}(I_i,T_{i,\text{s}})  P_{i}^{\text{ele}}(I_i,T_{i,\text{s}}) }{2 F U_{i}^{\text{cell}} (I_i,T_{i,\text{s}}) } }, \label{eq:nh2relax} \\
	\dot{n}_{i}^{\text{O}_2,\text{prod}} (I_i,T_{i,\text{s}}) \le {\color{black} \dfrac{  \eta_{i}^{\text{cell}}(I_i,T_{i,\text{s}})  P_{i}^{\text{ele}}(I_i,T_{i,\text{s}}) }{4 F U_{i}^{\text{cell}} (I_i,T_{i,\text{s}}) }}, \label{eq:no2relax}
\end{gather}
\noindent
then it is easy to find that the feasible space defined by (\ref{eq:nh2relax})--(\ref{eq:no2relax}) is convex. Further considering the control objective (\ref{eq:obj}), that is, to maximize $\dot{n}_{i}^{\text{H}_2,\text{prod}}$, and $\dot{n}_{i}^{\text{O}_2,\text{prod}}$, when the NMPC obtains the optimal solution, (\ref{eq:nh2relax})--(\ref{eq:no2relax}) all take the equality, that is, the above relaxation is exact.

\begin{figure}[t]
	\centering
	\includegraphics[scale=0.90]{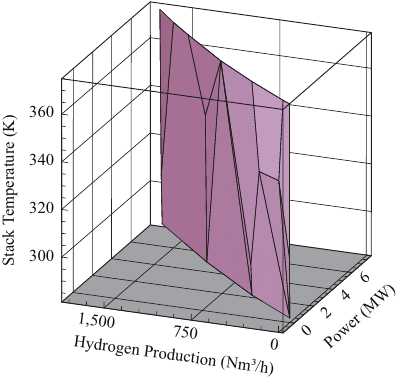}\vspace{-6pt}
	\caption{{\color{black}Polyhedral approximation of the production function of the 1,000 Nm$^3$/h-rated stack.}}
	\label{fig:dd}
\end{figure}

Then, we use the double description (DD) algorithm \cite{jones2010polytopic} to find the optimal polyhedral approximation of the feasible space defined by (\ref{eq:nh2relax})--(\ref{eq:no2relax}), formulated as a linear inequality. Due to $\dot{n}_{i}^{\text{H}_2,\text{prod}} = 2 \dot{n}_{i}^{\text{O}_2,\text{prod}}$, we only need to address the hydrogen output feasibility space, as:
\begin{gather}\color{black}
	 \bm{A} \Big[ P_i^{\text{ele}}, T_{i,\text{s}} , \dot{n}_{i}^{\text{H}_2,\text{prod}} \Big]^{\text{T}} + \bm{b} \le \bm{0}, \label{eq:prodieq}
\end{gather}
where $\bm{A}$ and $\bm{b}$ are constant 
matrix and vector.
Using a 1,000 Nm$^3$/h-rated stack studied in Section \ref{sec:case} as an example, with an error tolerance taken as 0.1\%, Fig. \ref{fig:dd} shows the polyhedral approximation. 
We use (\ref{eq:prodieq}) to replace the original model (\ref{eq:UItotal})--(\ref{eq:power}) in the controller, thereby eliminating non-convexity. 

\subsubsection{Reduced-Order Approximation of the HTO Accumulation Process}
\label{sec:htosimply}

As discussed in Section \ref{sec:hto}, the process of HTO accumulation takes three stages.
The time constants for the anode half-cell and liquid phase of the separator are on the order of seconds and minutes, while the gas-phase time constant ranges to tens of minutes. Thus, we can focus solely on the gas phase and adopt a first-order model:
\begin{align}
	\dfrac{\text{d} n^{\text{H}_2,\text{sep}}_{\text{gas}}}{\text{d}t} = \sum\nolimits_{i=1}^N \dot{n}^{\text{H}_2,\text{im}}_{i} - \dot{n}^{\text{H}_2,\text{im}}_{\text{out}}. \label{eq:htosimple}
\end{align}

Replacing (\ref{eq:Nbalance}) with (\ref{eq:htosimple}), the order of the model is reduced from $N$+2 to 1, improving both computational efficiency and numerical stability. Additionally, the molar quantity of HTO impurities in the half-cells and liquid phase, $ n^{\text{H}_2,\text{an}}_i$ and $n^{\text{H}_2,\text{sep}}_\text{liq}$, cannot be measured. After simplification, these variables are omitted, with the remaining $n^{\text{H}_2,\text{sep}}_\text{gas}$ easily measurable, thereby eliminating the need for a state estimator.

\begin{remark}
	The time constants of the heat transfer processes among stacks, separators, and the heat exchanger do not differ by order of magnitude. Therefore, unlike the HTO process, reducing the order of the temperature dynamic model would introduce noticeable errors. This has been discussed in studies on 1-in-1 systems \cite{qi2023thermal}, and the same applies to $N$-in-1 systems.
\end{remark}

\subsubsection{Bilinear Terms in Temperature and Mass Transfer Dynamics}
\label{sec:bilinear}

The state-space model presented in Section \ref{sec:model} has bilinear terms, which cannot be relaxed as convex as in Section \ref{sec:polynomial}. To address these terms, we employ a discretization-based big-M method \cite{bemporad1999control}.

The bilinear terms requiring reformulation include $I_i T_{i,\text{s}}$, $I^2_i$, $I_i {n}_{\text{gas}}^{\text{H}_2,\text{sep}}$, $v_{i,\text{lye}} T_{i,\text{s}}$, $v_{i,\text{lye}} T_{\text{sep,out}}$, and $v_{\text{c}} T_{\text{c,out}}$. Taking $I_i T_{i,\text{s}}$ as an example, the following steps are applied. First, discretize $I_i$ as
\begin{align}
  I_i = \sum\nolimits_{k=1}^{N^\text{d}} 2^{k-1} \beta^{I}_{i,k} \Delta I,   \label{eq:idist}
\end{align}
\noindent
where $N^\text{d}$ is the number of binary bits; $\Delta I = \overline{I}/2^{N^\text{d}}$ is the step size; and $\beta^{I}_{i,k}$ represents the $k$th binary variable. Then, $I_i T_{i,\text{s}}$ can be replaced with:
\begin{gather}
  I_i T_{i,\text{s}} = \sum\nolimits_{k=1}^{N^\text{d}} 2^{k-1} \delta^{I,T}_{i,k} \Delta I,   \label{eq:itdist} \\
  T_{i,\text{s}} - M \big(1- \beta^{I}_{i,k} \big) \le  \delta^{I,T}_{i,k}  \le T_{i,\text{s}} + M  \big(1- \beta^{I}_{i,k} \big), \label{eq:itbigm1} \\
  - M  \beta^{I}_{i,k}  \le  \delta^{I,T}_{i,k}  \le M  \beta^{I}_{i,k} , \label{eq:itbigm2}
\end{gather}
where $\delta^{I,T}_{i,k}$ is an intermediate variable; $M$ is a sufficiently large constant. This reformulation converts bilinear terms into mixed-integer linear constraints, which can be efficiently solved using commercial solvers.

\subsection{Controller Summary}
\label{sec:controlsummary}

After the aforementioned simplifications, {\color{black}the original NMPC is reformulated and implemented as an MIQP-based MPC problem}, expressed compactly as:
\begin{align}
	{\color{black}\min}\ (\ref{eq:obj}), \ \text{s.t.} \ (\ref{eq:diststate}), (\ref{eq:prodieq}), (\ref{eq:limitdcpower}){\text{--}}(\ref{eq:limitcool}),(\ref{eq:idist}){\text{--}}(\ref{eq:itbigm2}).    \label{eq:mpcall}
\end{align}

With this formulation, we can set a control horizon ranging from a few minutes to an hour, and the step length to several minutes. Each optimization can be solved within seconds.

In implementation, (\ref{eq:mpcall}) is solved repeatedly. At each step, only the first control action is applied, and the process is updated with new measurements. By repeating this process, the controller distributes power among stacks and adjusts lye and cooling flows in real time, compensating for disturbances and model errors.
The pressure and liquid levels are regulated by a separate controller, which is not covered in this work.

\begin{remark}
	\label{remark:2}
	The proposed {\color{black}MIQP-MPC} applies to both $N$-in-1 and 1-in-1 systems and can be adapted to different lye circulation topologies. This provides a consistent basis for performance comparison across different configurations, as discussed in Sections \ref{sec:comparetopo}--\ref{sec:comparemulti}. Note that the focus of this work is not to benchmark or optimize control performance against conventional control methods.
\end{remark}

\section{Simulation Studies}
\label{sec:case}

\subsection{Simulation System Setup}
\label{sec:setting}

An experiment directly implementing the proposed {\color{black}MIQP-MPC} on the industrial-scale system is infeasible due to  potential safety risks, production interruption, and costs. Therefore, the calibrated 4,000 Nm$^3$/h-rated 4-in-1 system model from Section \ref{sec:validation} (parameters in Table \ref{tab:parameter}, Appendix B) is used to evaluate control performance. To compare lye circulation designs, three models with 4-pump, 2-pump, and 1-pump topologies are constructed using the same parameters.

The controller uses a horizon $N^\text{h}=1,800$ s and a step length $\Delta t =  450$ s. The weight coefficients in the objective (\ref{eq:obj}) follow $\lambda^{\text{prod}}=1$,  $ \lambda^{\text{track}}=1.2$, $\lambda^{\text{temp}}=0.15$, $\lambda^{\text{I}}=0.0002$, $ \lambda^{\text{lye}}  = 25,000$, and $ \lambda^{\text{c}} = 0.5$. The complete nonlinear model presented in Table \ref{tab:model} is used for time-domain simulation. The controller updates the control commands every 10 seconds.
Simulations are performed on \emph{Wolfram Mathematica 13.0}, and the optimization problem in the controller is solved via \emph{Gurobi 11.0.3}. {\color{black}The computational platform is a desktop equipped with an \emph{Intel Core i7-12700K}  CPU and a RAM of 32.0 GB. }

It should be noted that the objective of this study is to assess the relative dynamic performance of 1-in-1 and $N$-in-1 configurations rather than benchmarking {\color{black}MIQP-MPC} against conventional control strategies. The {\color{black}MIQP-MPC} is applied uniformly across both configurations to enable fair comparison; the observation in dynamic performance hence reflects system characteristics rather than the superiority of the control strategy itself.

\begin{figure}[tb]
	\centering
	\includegraphics[scale=0.94]{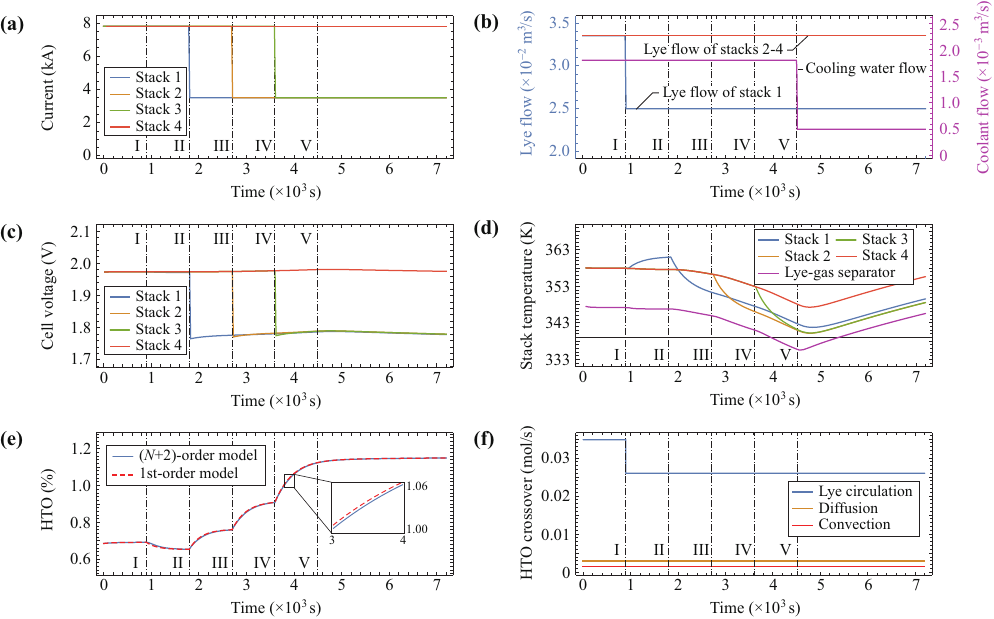}\vspace{-6pt}
	\caption{Open-loop responses of the 4-in-1 AWE system. (a) Electrolytic current. (b) Lye and cooling water flow rates. (c) Cell voltage. (d) Stack and separator temperatures. {\color{black}(e) HTO impurity at the gas phase of the separator, obtained by the ($N$+2)th-order model and the $1$st-order model. (f) Components of HTO crossover in Stack 1.}}
	\label{fig:openloop}
\end{figure}

\subsection{Open-Loop Simulation}
\label{sec:openloop}

An open-loop simulation (without feedback control) illustrates the coupling between inputs and system responses shown in Fig. \ref{fig:state}. The 4-pump topology in Fig. \ref{fig:cir}(a) is used, allowing independent lye flow control. Initial conditions are set at stack temperature 85 $^\circ$C (358 K), separator temperature 74 $^\circ$C (346 K), and HTO concentration 0.69\%.
The control actions for current, lye flow, and cooling water flow are plotted in Figs. \ref{fig:openloop}(a) and \ref{fig:openloop}(b).

At 900 s, when the lye flow rate of stack 1 decreases from $3.35 \times 10^{-2}$ m$^3$/s to $2.5 \times 10^{-2}$ m$^3$/s (Action I), its temperature rises by 0.5 K, while HTO slightly decreases. At 1,800 s, 2,700 s, and 3,600 s, currents of stacks 1--3 drop from 7,800 A to 3,500 A (Actions II--IV). Stack temperatures decrease, the separator temperature falls to 335.7 K (62.7 $^\circ$C), and HTO increases from 0.66\% to 1.15\%.
At 4,500 s, cooling water flow decreases from $1.78 \times 10^{-3}$ m$^3$/s to $0.50 \times 10^{-3}$ m$^3$/s (Action V), raising all temperatures by about 11 K within 2,500 s. {\color{black}During the whole process, lye circulation contributes more than 80\% of the total HTO crossover, as shown in Fig. \ref{fig:openloop}(f), which is consistent with existing findings \cite{qi2021pressure, guan2026region}.}
These results reveal strong interactions among variables, confirming the need for coordinated MIMO control.

{\color{black}
Fig. \ref{fig:openloop}(e) further compares the full-order HTO impurity model presented in Section \ref{sec:hto} with the simplified first-order model introduced in Section \ref{sec:htosimply}. The two models exhibit very similar responses, while the first-order model responds slightly faster during transient conditions. However, the difference remains very small, with the maximum mismatch being less than 0.01\%. Therefore, the reduced-order model can be adopted in the controller formulation while maintaining satisfactory safety and accuracy.
}

\begin{figure}[tb]
	\centering
	\includegraphics[scale=0.94]{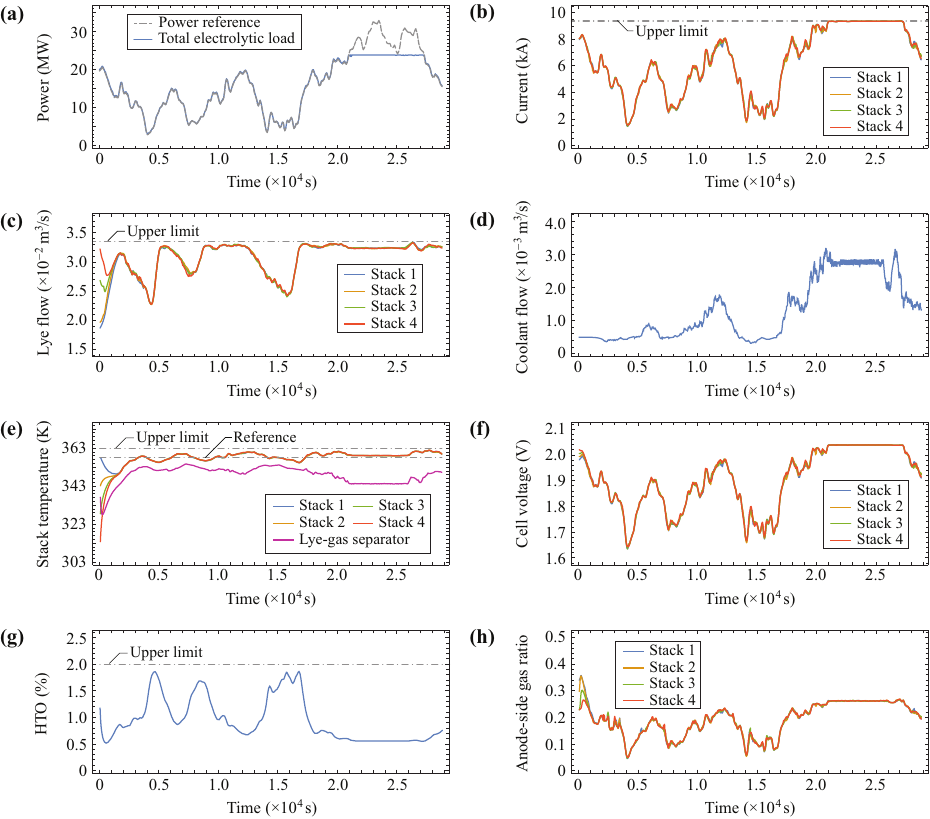}\vspace{-6pt}
	\caption{{\color{black}Control and responses of the 4-in-1 AWE system with a 4-pump configuration. (a) Load power reference and power consumption. (b) Electrolytic currents. (c) Lye flow rates. (d) Cooling water flow rate. (e) Stack and separator temperature. (f) Cell voltages. (g) HTO impurity. (h) Gas ratio in the anode half-cells.}}
	\label{fig:4valve}
\end{figure}

\subsection{Overall Performance under Dynamic Operation}
\label{sec:casebase}

The performance of the 4-in-1 system under the proposed {\color{black}MIQP-MPC}
is analyzed. The simulation adopts the 4-pump configuration shown in Fig. \ref{fig:cir}(a).
The power reference $P_{\text{tot}}^{\text{ref}}(t)$ emulates the fluctuations of renewable energy and varies between 6 MW and 38 MW over 8 hours, as shown in Fig. \ref{fig:4valve}(a), with low power in the first 5 hours and higher levels thereafter.

The initial temperatures of stacks are set to 85 $^\circ$C, 70 $^\circ$C, 55 $^\circ$C, and 40 $^\circ$C (358 K, 343 K, 328 K, and 313 K), respectively. The separator temperature is set to 65 $^\circ$C (338 K), with an initial HTO concentration at 1.2\%. Throughout the simulation, the controller continuously adjusts the current, lye flow, and cooling water flow for each stack, as shown in Figs. \ref{fig:4valve}(b)--(d). {\color{black}The minimal, maximal, and average computation time for solving the MIQP-based control problem are 1.011, 7.296, and 2.769 s.} The behaviors of the 4-in-1 system are discussed below.

\emph{a) Inter-stack Load Allocation:} The load is evenly shared among the 4 stacks throughout the process, driven by the objective function (\ref{eq:obj}). As shown in Fig. \ref{fig:dd}, the relation between hydrogen yield, electrolytic power, and temperature is concave.  At optimality, marginal production rates are equal across stacks, i.e., $\partial \dot{n}_{i}^{\text{H}_2,\text{prod}} / \partial P_{i}^{\text{ele}} = \partial \dot{n}_{j}^{\text{H}_2,\text{prod}} / \partial P_{j}^{\text{ele}}$, $\forall i,j\le N$. With identical stacks, this yields uniform loading, known as the \emph{equimarginal principle} \cite{li2024two,zeng2024scheduling}. After 20,000 s, when total power exceeds 24 MW, all stacks operate at full load.

\emph{b) Temperature Dynamics:} Due to coupling through the BoP, stack temperatures converge to the setpoint 85 $^\circ$C.
As shown in Figs. \ref{fig:4valve}(c) and \ref{fig:4valve}(e), stacks 4 and 3, initially cooler than the separator, increase their lye flow rates to accelerate heating, while stacks 1 and 2 reduce lye flow rates to minimize heat loss.
The minimum flow rates for stacks 1 and 2 reach 0.0186 m$^3$/s and 0.0197 m$^3$/s compared to the rated value of 0.0345 m$^3$/s. Heat exchange rapidly makes stack temperatures become uniform during startup.
After convergence, lye flow remains near rated values to stabilize gas-liquid ratios, and temperature control relies mainly on cooling water. Stack temperatures remain within [83, 88] $^\circ$C.

\emph{c) HTO Accumulation Dynamics:} Since HTO impurity is enforced via constraints rather than the objective, its fluctuations do not affect the control actions in most cases.
However, during significant load drops (e.g., at 4,000 s, 8,000 s, and 15,000 s), HTO rises rapidly, similar to 1-in-1 systems.
To keep it below 2\%, the controller reduces lye flow to 0.0226, 0.0285, and 0.0242 m$^3$/s (67\%, 85\%, and 72\% of rated), limiting cross-mixing as described by (\ref{eq:nlye}).
Otherwise, lye flow stays near nominal values, as shown in Fig. \ref{fig:4valve}(c).

Overall, despite strong coupling between production, temperature, and HTO dynamics in the $N$-in-1 system, the proposed controller can manage these complexities effectively, thereby answering Question 2 posed in Section \ref{sec:intro}.

\subsection{Comparisons Between Different Lye Circulation Topologies}
\label{sec:comparetopo}

\begin{figure}[tb]
	\centering
	\includegraphics[scale=0.94]{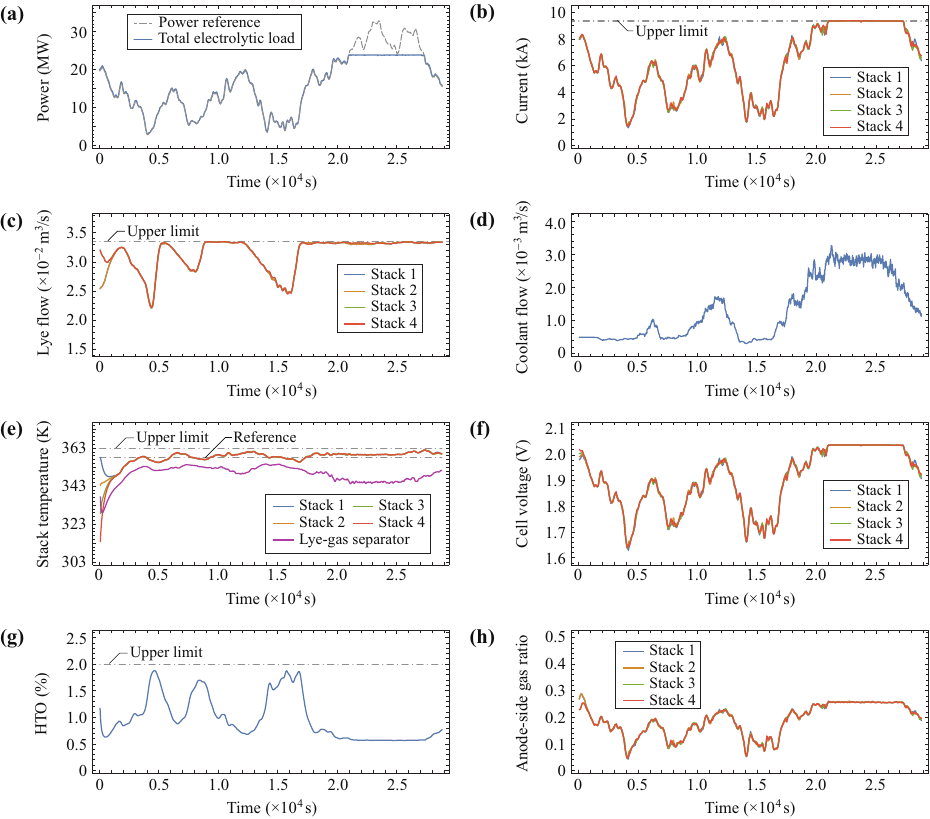}\vspace{-6pt}
	\caption{{\color{black}Control and responses of the 4-in-1 AWE system with a 2-pump configuration shown in Fig. \ref{fig:cir}(b). (a) Load power reference and power consumption. (b) Electrolytic currents. (c) Lye flow rates. (d) Cooling water flow rate. (e) Stack and separator temperatures. (f) Cell voltage of the stacks. (g) HTO impurity. (h) Gas ratio in the anode half-cells.}}
	\label{fig:2valve}
\end{figure}

\begin{figure}[tb]
	\centering
	\includegraphics[scale=0.94]{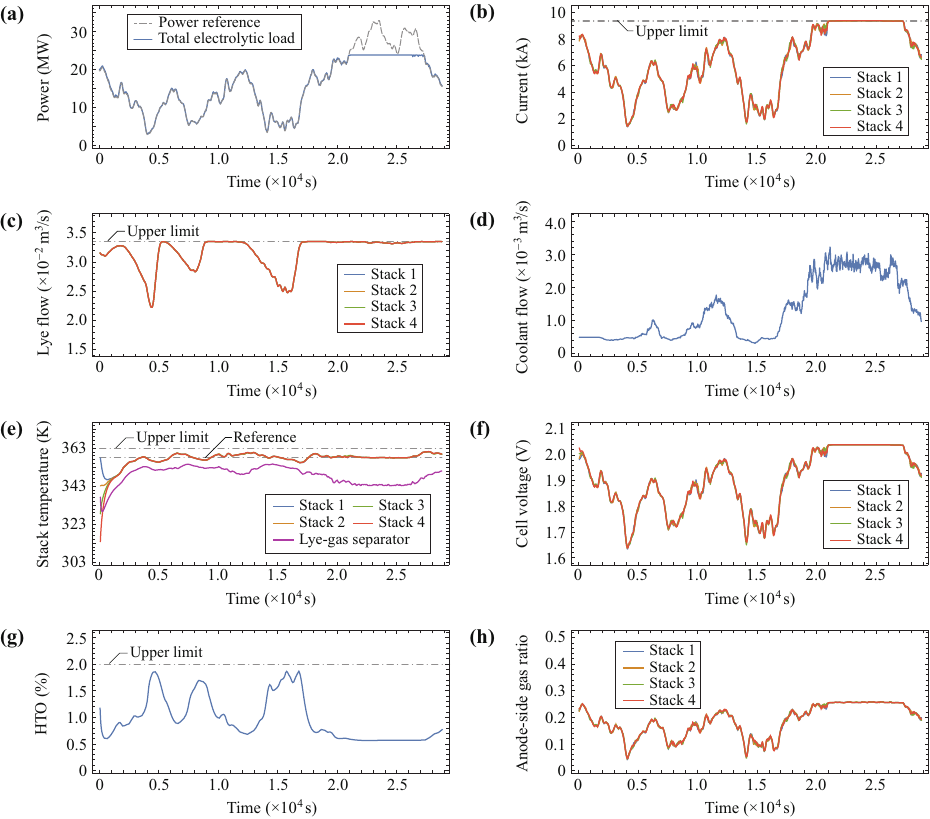}\vspace{-6pt}
	\caption{{\color{black}Control and responses of the 4-in-1 AWE system with a 1-pump configuration shown in Fig. \ref{fig:cir}(c). (a) Load power reference and power consumption. (b) Electrolytic currents. (c) Lye flow rates. (d) Cooling water flow rate. (e) Stack and separator temperature. (f) Cell voltage of the stacks. (g) HTO impurity. (h) Gas ratio in the anode half-cells.}}
	\label{fig:1valve}
\end{figure}

\begin{figure}[tb]
	\centering
	\includegraphics[scale=0.94]{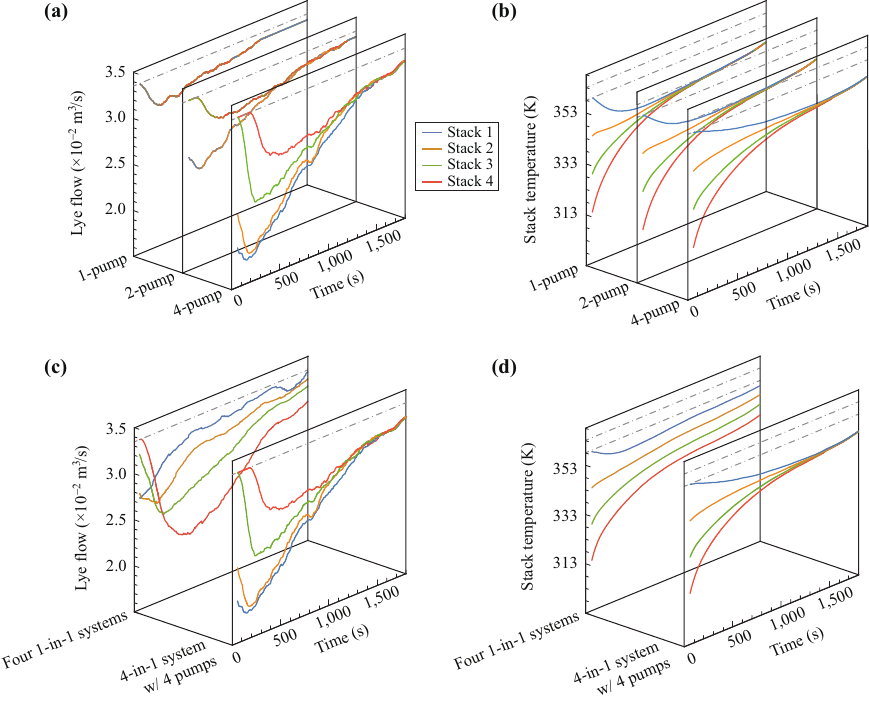}\vspace{-6pt}
	\caption{Comparison on control and responses between the 4-in-1 AWE system with  4-pump, 2-pump, and 1-pump configurations and four traditional 1-in-1 system in parallel. (a) and (c) Lye flow rates of the stacks. (b) and (d) Temperature of the stacks.}
	\label{fig:compcir}
\end{figure}

\begin{table}[tb]\footnotesize\centering
	\renewcommand{\arraystretch}{1.45}
	\caption{Performance comparison between 4-in-1 AWE systems with different lye circulation topologies and four 1-in-1 systems operating in parallel}\vspace{6pt}
	\label{tab:compare}
	\begin{tabular}{ cccccc }\hline\hline
		Configuration & \tabincell{c}{Energy use\\(MWh)}     & \color{black}\tabincell{c}{Load-tracking \\ RMSE for unsaturated \\  periods (MW)}                & \tabincell{c}{RMSE for stack \\   temperature  \\control (K)} &  \tabincell{c}{Total hydrogen \\ yield (Nm$^3$) }  &  \tabincell{c}{Specific energy \\consumption \\(kWh/Nm$^3$) } \\ \hline
		4-in-1 (4-pump)       & \color{black}123.657   	  & \color{black}0.1636	  & \color{black}3.757   & \color{black}25,847.6  & \color{black}4.784 \\
		4-in-1 (2-pump)       & \color{black}123.668	  & \color{black}0.1641	  & \color{black}3.814   & \color{black}25,848.8	& \color{black}4.784 \\
		4-in-1 (1-pump)       & \color{black}123.691	  & \color{black}0.1649	  & \color{black}3.765	& \color{black}25,851.6	& \color{black}4.785 \\ \hline
		\tabincell{c}{Four 1-in-1 systems\\{\color{black}(all online)}}   & \color{black}123.741	  & \color{black}0.1704	  & \color{black}4.177	& \color{black}25,862.3	& \color{black}4.78{\color{black}5} \\
		\color{black} \tabincell{c}{Four 1-in-1 systems\\(with on-off switching)}   & \color{black} 123.733	  & \color{black}0.1710	  & \color{black} 8.394	& \color{black} 25,475.3	&  \color{black} 4.857  \\ \hline\hline
	\end{tabular}
\end{table}

As discussed in Section \ref{sec:lye}, the degree of freedom in controlling lye flow distribution varies across different topologies.
Their impact is evaluated using the three configurations in Fig. \ref{fig:cir}, with identical system parameters and power reference shown in Fig. \ref{fig:4valve}(a).
Figs. \ref{fig:2valve} and \ref{fig:1valve} show results for the 2-pump and 1-pump systems. Apart from differences in inter-stack lye flow distribution during startup, current control, cooling flow, and performance metrics (power tracking, temperature, and HTO) are nearly identical across all topologies.

To further analyze the transient behavior of different topologies, the control of lye flow and the temperature responses during the first 30 minutes are shown in Figs. \ref{fig:compcir}(a) and \ref{fig:compcir}(b).
In the 4-pump topology, lye flow rates in the stacks decrease by large and differing levels first, as analyzed in Section \ref{sec:casebase}.
In contrast, in the 2-pump and 1-pump systems, stacks sharing a common circulation pump exhibit nearly identical lye flow rates, which cannot be independently adjusted, as noted in Section \ref{sec:lye} and Remark \ref{remark:1}.
Therefore, these systems cannot actively regulate inter-stack heat exchange through differential lye flow, and the controller instead maintains higher lye flow rates (the minimum to be 0.0216 m$^3$/s and 0.0305 m$^3$/s for the 2-pump and 1-pump configurations, respectively) across all stacks to accelerate the process of heating up.

Consequently, temperature drops of stacks 1 and 2 in the 2-pump and 1-pump configurations  (10.5 K and 12 K for stack 1 in the 2-pump and 1-pump systems) are larger than the 4-pump system (only 8.5 K due to independent lye flow control).
After startup, all temperatures converge to 85 $^\circ$C and flow rates return to nominal values, with no significant differences thereafter, as shown in Figs. \ref{fig:4valve}, \ref{fig:2valve}, and \ref{fig:1valve}.

Table \ref{tab:compare} summarizes the key performance metrics. {\color{black}The load-tracking RMSEs are calculated for the first 5 hours because in the remaining 3 hours the tracking error is dominated by saturation, as shown in Fig. \ref{fig:4valve}--\ref{fig:1valve}.} Despite the reduced control flexibility in the 2-pump and 1-pump configurations, the differences in specific energy consumption and temperature control error are less than 0.001 kWh/Nm$^3$ and 0.05 K. {\color{black}The load-tracking RMSEs during unsaturated periods also remain comparable, with a maximum difference of only 0.013 MW.}
These results show that while topology affects transient behavior, the proposed controller maintains comparable performance under the continuous-operation scenarios considered here.

{\color{black}
Note that the performance differences among the compared topologies are obtained under the same calibrated simulation framework and mainly reflect relative performance differences between topologies. These differences should be interpreted as comparative trends rather than exact absolute temperature prediction improvements relative to industrial measurements.
}

\subsection{Comparing the $N$-in-1 System with Multiple 1-in-1 AWE Systems}
\label{sec:compare1in1}

\begin{figure}[tb]
	\centering
	\includegraphics[scale=0.94]{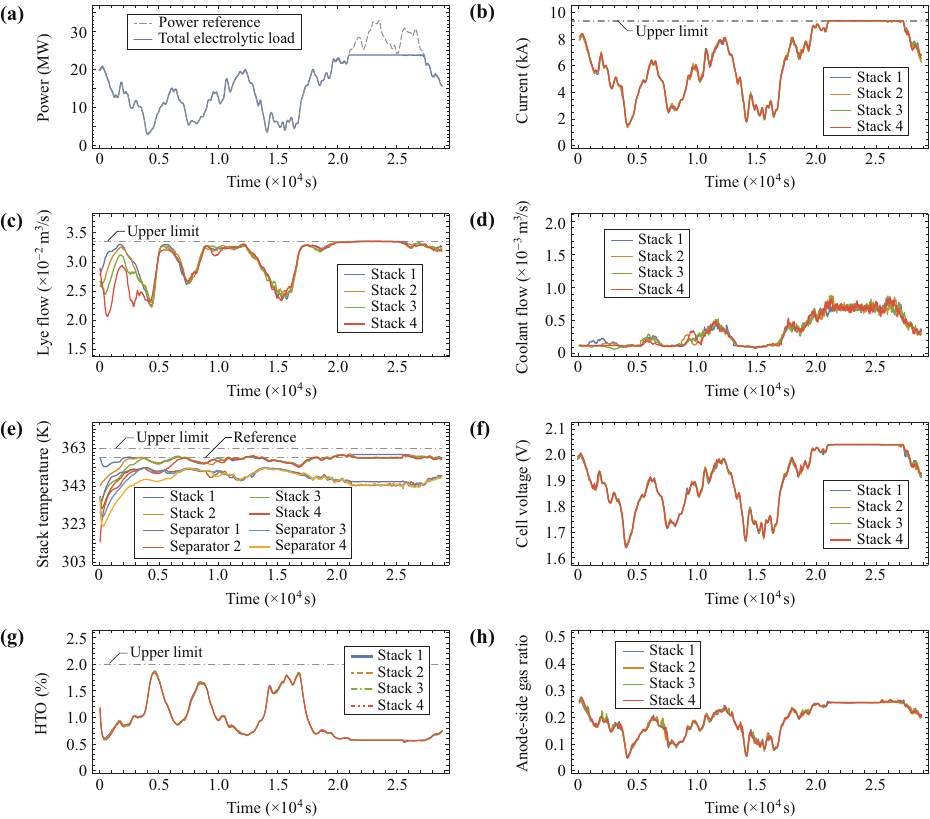}\vspace{-6pt}
	\caption{Control and responses of four traditional 1-in-1 AWE systems operating in parallel {\color{black} and all remaining online}. (a) Load power reference and power consumption. (b) Electrolytic currents. (c) Lye flow rates. (d) Cooling water flow rate. (e) Stack and separator temperatures. (f) Cell voltage of the stacks. (g) HTO impurities. (h) Gas ratios in the anode half-cells.}
	\label{fig:four1in1}
\end{figure}

\begin{figure}[tb]
	\centering
	\includegraphics[scale=0.94]{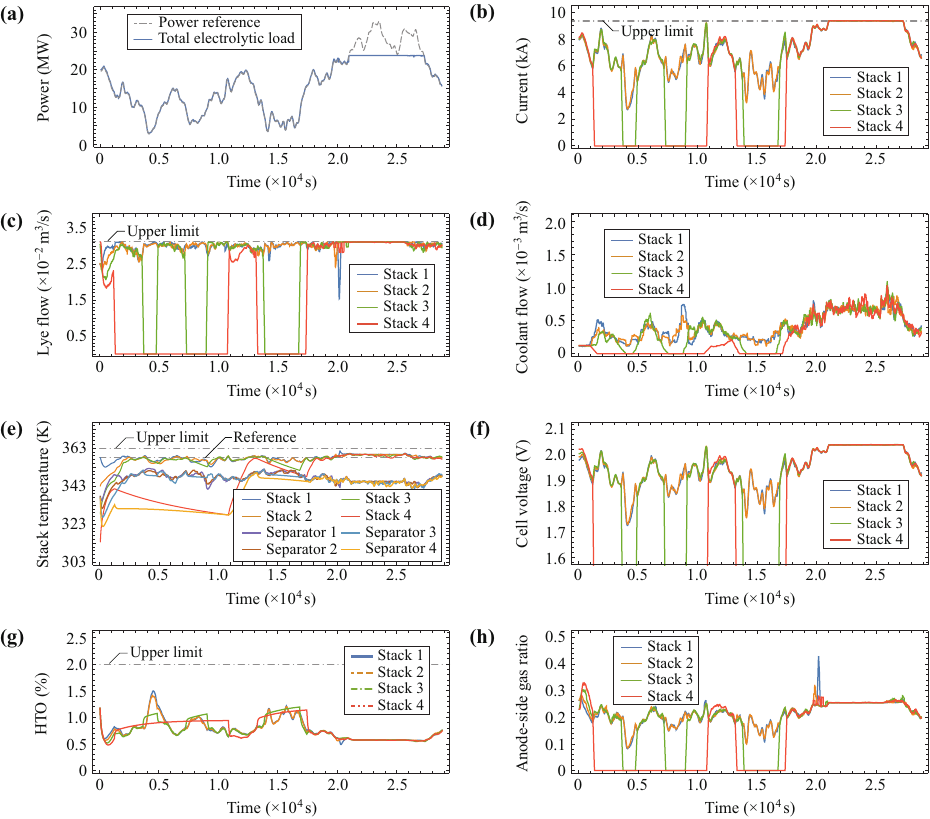}\vspace{-6pt}
	\caption{{\color{black}Control and responses of four traditional 1-in-1 AWE systems operating in parallel with switching between online and offline states. (a) Load power reference and power consumption. (b) Electrolytic currents. (c) Lye flow rates. (d) Cooling water flow rate. (e) Stack and separator temperatures. (f) Cell voltage of the stacks. (g) HTO impurities. (h) Gas ratios in the anode half-cells.}}
	\label{fig:four1in1onoff}
\end{figure}

To address Question 1 in Section \ref{sec:intro}, we compare $N$-in-1 and traditional 1-in-1 AWE systems in terms of flexibility and efficiency. Although many ReP2H projects adopt $N$-in-1 configurations, most studies focus on 1-in-1 systems. A systematic comparison is therefore needed.

We compare the 4,000 Nm$^3$/h-rated 4-in-1 system  in Section \ref{sec:setting} with four parallel 1,000 Nm$^3$/h-rated 1-in-1 systems under the same load reference in Fig. \ref{fig:4valve}(a). The 1-in-1 parameters are based on the same experimental system, with identical stack settings and proportionally scaled BoP components. For consistency, the same {\color{black}MIQP-MPC} is applied with $N=1$ and unchanged settings.

As in Section \ref{sec:casebase}, the initial temperatures of the stacks in the 1-in-1 systems are set to 85~$^\circ$C, 70~$^\circ$C, 55~$^\circ$C, and 40~$^\circ$C, respectively, with the initial temperature of the separators to be 65 $^\circ$C, and the initial HTO concentration at 1.2\%.

  {\color{black}
    Since 1-in-1 systems can be started up and shut down independently, two operating settings are considered for the four independent 1-in-1 systems: a) all four 1-in-1 systems remain online; and b) the 1-in-1 systems are allowed to switch between online and offline states every 10 minutes according to the available power. The number of online systems is determined by $ N^\text{on} = \max\{ 0 , \min\{4, \lceil P_{\text{tot}}^{\text{ref}} / \overline{P}^\text{ele} \rceil \} \}$.

    Figs. \ref{fig:four1in1} and  \ref{fig:four1in1onoff} present the simulation results under settings a) and b), respectively. Figs. \ref{fig:compcir}(c) and \ref{fig:compcir}(d) compare the lye flow control and temperature responses during the initial phase before on-off switching occurs.

    For the all-online setting, the 4-in-1 system and the four 1-in-1 systems show similar behavior over the 8-hour simulation, except during startup. As discussed in Section \ref{sec:casebase}, the equimarginal principle leads to nearly uniform power sharing to maximize hydrogen production. For the four 1-in-1 systems, independent temperature and lye-flow control results in similar control actions under the same objective.

    When on-off switching is allowed, the load-tracking performance remains comparable to that of the all-online setting. Since two systems are switched off during low-load periods, the online systems operate at a higher average load, and the HTO impurity concentration becomes lower, indicating a larger safety margin. However, this benefit comes with additional thermal transients. The offline stacks cool down during low-load periods, which increases the temperature-tracking error and reduces energy efficiency.
    }

    Another difference arises during startup. In the 4-in-1 system, inter-stack heat exchange via differential lye flow accelerates temperature equalization. In contrast, the 1-in-1 systems lack such coupling, leading to slower convergence. In the 4-in-1 case, the largest lye flow reduction occurs in the hottest stack (up to 44.5\%), while in the 1-in-1 systems it occurs in the coldest stack (37.3\%). After startup, temperatures converge in both cases, but more slowly in the 1-in-1 systems due to thermal isolation.

    Table \ref{tab:compare} summarizes the performance metrics.
    {\color{black}
    For setting a), where all four 1-in-1 systems remain online, the 4-in-1 system and the four 1-in-1 systems show comparable efficiency, load-tracking performance, and temperature regulation. Differences in temperature RMSE and specific energy consumption are below 0.420 K and 0.001 kWh/Nm$^3$, respectively, and the load-tracking RMSE during unsaturated periods differs by only 0.005 MW. This similarity reflects the same optimal allocation driven by electrochemical characteristics summarized in Section \ref{sec:production}. Therefore, in response to Question 1 under continuous all-online operation, adopting the 4-in-1 system does not reduce operational flexibility or energy efficiency.
    }

    {\color{black}
    For setting b), where independent 1-in-1 systems are allowed to switch between online and offline states, load-tracking performance remains comparable. However, the temperature-tracking RMSE increases to 8.394 K, the specific energy consumption increases from 4.785 to 4.857 kWh/Nm$^3$, and the total hydrogen yield decreases from 25,862.3 to 25,475.3 Nm$^3$. These results suggest that modular on-off operation can improve low-load flexibility and reduce HTO impurity, but may also introduce additional thermal transients and efficiency penalties if the startup/shutdown schedule is not optimized.

    It should be noted that the on-off scheduling of multiple AWE systems is a complex hybrid optimization problem \cite{varela2021modeling,qiu2023extended,zeng2024scheduling,li2024two}. A multi-layer energy management framework is needed to coordinate electrical dispatch, startup/shutdown decisions, pressure and liquid-level constraints, temperature dynamics, and HTO impurity limits in large hydrogen plants with multiple 1-in-1 or $N$-in-1 systems. This remains an important challenge for addressing Question 2 under highly volatile operating conditions and merits future study.
    }

\subsection{Comprehensive Comparisons under Various Scenarios}
\label{sec:comparemulti}

\begin{figure}[tb]
	\centering
	\includegraphics[scale=0.94]{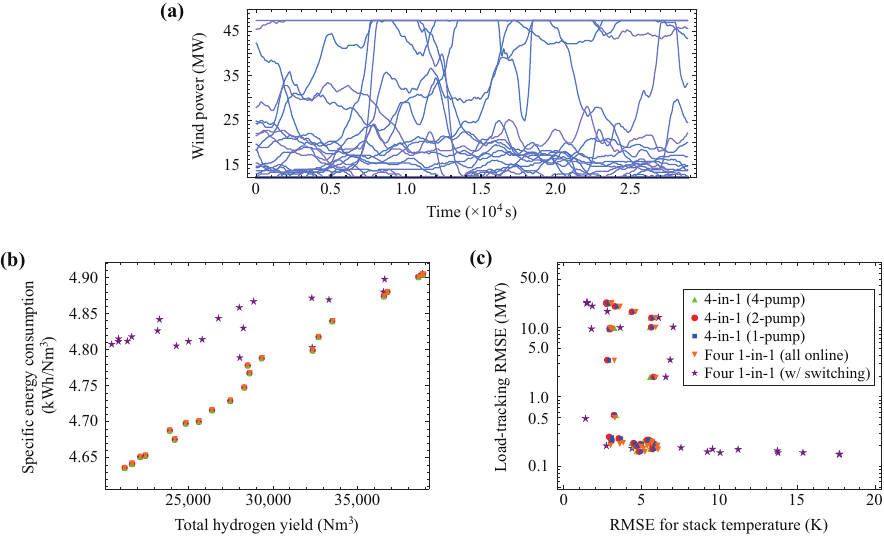}\vspace{-6pt}
	\caption{Performance metrics of the 4-in-1 AWE systems with different lye circulation topologies and four 1-in-1 systems operating in parallel under 25 scenarios of wind power supply. (a) Power reference signals based on a wind farm by the Ris{\o}e. (b) Total hydrogen yield and specific energy consumption. (c) RMSE for stack temperature and load-tracking RMSE.}
	\label{fig:multi}
\end{figure}

To further compare $N$-in-1 and 1-in-1 systems, we simulate 25 wind power scenarios. Each scenario (Fig. \ref{fig:multi}(a)) is derived from offshore wind data reported by Ris{\o} \cite{lin2012assessment}. Performance metrics, including total energy use, load-tracking RMSE, temperature RMSE, hydrogen yield, and specific energy consumption, are shown in Figs. \ref{fig:multi}(b)--\ref{fig:multi}(c), with averages in Table \ref{tab:multi}.

\begin{table}[tb]\footnotesize\centering
	\renewcommand{\arraystretch}{1.4}
	\caption{{\color{black}Average performance metrics of the 4-in-1 AWE systems with different lye circulation topologies and four 1-in-1 systems operating in parallel under 25 scenarios of wind power supply}}\vspace{6pt}
	\label{tab:multi}
	\begin{tabular}{ cccccc }\hline\hline
	Configuration & \tabincell{c}{Energy use\\(MWh)}        & \color{black}\tabincell{c}{Load-tracking \\ RMSE for unsaturated \\  periods (MW)}   & \tabincell{c}{RMSE for stack \\   temperature \\ control (K)} &  \tabincell{c}{Total hydrogen \\ yield (Nm$^3$) }  &  \tabincell{c}{Specific energy \\consumption \\ (kWh/Nm$^3$) } \\ \hline
	4-in-1 (4-pump)       & \color{black}142.214  & \color{black}0.2043  & \color{black}4.270   & \color{black}29,688.2  & \color{black}4.790 \\
	4-in-1 (2-pump)       & \color{black}142.208	 & \color{black}0.2088 & \color{black}4.328   & \color{black}29,686.4	& \color{black}4.790 \\
	4-in-1 (1-pump)       & \color{black}142.223	 & \color{black}0.2102 & \color{black}4.425	& \color{black}29,688.9	& \color{black}4.790 \\ \hline
	\tabincell{c}{Four 1-in-1 systems\\{\color{black}(all online)}}   & \color{black}142.282	& \color{black}0.1929 	& \color{black}4.614	& \color{black}29,697.7	& \color{black}4.791 \\
	\color{black} \tabincell{c}{Four 1-in-1 systems\\(with on-off switching)}   &  \color{black}142.121		& \color{black}0.1542  &\color{black} 7.068	  &   \color{black}29,270.5	& 	\color{black}4.855    \\ \hline\hline
\end{tabular}
\end{table}

The results confirm the findings of Sections \ref{sec:comparetopo} and \ref{sec:compare1in1}. Despite reduced degrees of freedom in lye flow control, the 4-in-1 system maintains comparable performance through coordinated adjustment of current, lye flow, and cooling. Moreover, because the stacks in the 4-in-1 system share a large heat capacity, the temperature fluctuations can be better stabilized under varying power supply. {\color{black}Different from the 1-in-1 systems with on-off switching that exhibit deviations in efficiency and temperature control,} the average differences between the 4-in-1 systems and the {\color{black}all-online} 1-in-1 systems in stack temperature RMSE, and specific energy consumption are below {\color{black}$0.344$} K, and $0.001$ kWh/Nm$^3$, relatively small for the system-level comparison considered here.  In summary, we answer Question 1 by confirming that selecting the 4-in-1 system does not compromise flexibility when compared to traditional 1-in-1 systems {\color{black} in continuous production}.

{\color{black}
	
\subsection{Performance of Different Lye Circulation Topologies under Asymmetric Operation}
\label{sec:asymmetry}

In practical operation, one or more stacks in an $N$-in-1 AWE system may have to operate at a reduced power level due to degradation or maintenance requirements. Asymmetric operation may affect lye-flow allocation, temperature regulation, HTO crossover, and energy efficiency, especially when several stacks share a common circulation pump.

To investigate this issue, an additional case is considered. The current upper limit of Stack 4 is reduced to 0.6 times the rated current, while the other stacks remain available at their full rated current limits. The same power-reference profile as in Section \ref{sec:casebase}  is used. The 4-pump, 2-pump, and 1-pump configurations are compared, and the performance metrics over the 8-hour simulation are summarized in Table \ref{tab:asymmetry}.

The results show that the load-tracking performance is mainly affected by the reduced total available capacity. During unsaturated periods, the controller can still track the power reference with comparable accuracy. When the power reference exceeds the reduced available capacity, the tracking error is mainly caused by saturation.

The three topologies also show comparable temperature-regulation and efficiency performance under this derated-stack condition. Through stack-current redistribution and lye-flow adjustment, the MIQP-based controller can partly compensate for the asymmetric capacity constraint and the resulting thermal imbalance. Therefore, under the moderate single-stack derating level considered here, the performance comparison among the 4-pump, 2-pump, and 1-pump configurations is not qualitatively changed.

It should be noted that this case only represents a moderate abnormal operating condition. More severe derating, stack idling, or multiple unavailable stacks may lead to stronger topology-dependent differences. Handling these conditions involves coupled load allocation, startup/shutdown decisions, temperature, pressure and liquid level regulation. A comprehensive energy management framework for such conditions is beyond the scope of this work and will be investigated in future studies.

\begin{table}[tb]\footnotesize\centering
	\renewcommand{\arraystretch}{1.4}
	\caption{\color{black}Performance metrics of the 4-in-1 AWE systems with different lye circulation topologies when Stack 4 is derated }\vspace{6pt}
	\label{tab:asymmetry}
	\begin{tabular}{ cccccc }\hline\hline
		\color{black}Configuration &  \color{black}\tabincell{c}{Energy use\\(MWh)}        & \color{black}\tabincell{c}{Load-tracking \\ RMSE for unsaturated \\  periods (MW)}   &  \color{black}\tabincell{c}{RMSE for stack \\   temperature \\ control (K)} &   \color{black}\tabincell{c}{Total hydrogen \\ yield (Nm$^3$) }  &  \color{black} \tabincell{c}{Specific energy \\consumption \\ (kWh/Nm$^3$) } \\ \hline
		\color{black}4-in-1 (4-pump)       & \color{black} 118.443  & \color{black}0.1641  &  \color{black} 4.277   &  \color{black}24,830.3  & \color{black} 4.7701\\
		\color{black}4-in-1 (2-pump)       & \color{black}118.517	 & \color{black}0.1690&  \color{black}4.526  &  \color{black}24,843.0	& \color{black}4.7706 \\
		\color{black}4-in-1 (1-pump)       & \color{black}118.509	 & \color{black}0.1676 &  \color{black}4.233	&  \color{black}24,840.3	& \color{black}4.7707 \\ \hline\hline
	\end{tabular}
\end{table}

}

\section{Conclusions}
\label{sec:conclusion}

This paper develops a state-space model for $N$-in-1 AWE systems, i.e., multi-stack AWE systems with shared separation and lye circulation. The model captures coupled dynamics of electrochemical reactions, lye flow, temperature, and HTO accumulation. Experiments on an industrial 4,000 Nm$^3$/h-rated 4-in-1 system validate its accuracy.

Then, {\color{black}an NMPC framework and an MIQP-based solution method} are designed to coordinate current distribution, lye circulation, and cooling under fluctuating power input, primarily as a tool allowing for a fair comparison among different topology designs and 1-in-1 systems.
While experimental implementation of the {\color{black}MIQP-MPC} on the industrial 4-in-1 system is not feasible, the simulated control results are sufficient to compare the dynamic performance and flexibility between 1-in-1 and N-in-1 configurations.  Comprehensive simulation studies lead to the following findings:
\begin{enumerate}
	\item With appropriate control, $N$-in-1 systems with different lye circulation topologies achieve similar performance in efficiency, load tracking, and temperature regulation {\color{black}in both normal and partially derated conditions}, despite different degrees of freedom of control.
	\item  $N$-in-1 and multiple 1-in-1 systems exhibit similar performance {\color{black}under continuous production conditions.} Differences in efficiency, load tracking, and temperature control are relatively small for practical applications {\color{black}when all stacks are in production}.
	\item As in multiple parallel $1$-in-1 AWE systems, the stacks in the $N$-in-1 systems tend to share a common operational state, i.e., electrolytic power, temperature, and lye flow, to maximize hydrogen production driven by the equimarginal principle.
\end{enumerate}

Future research may address the following aspects. First, it is expected to explore scaled pilot tests or incremental implementation to further validate control strategies on industrial systems. Furthermore, this study assumes continuous operation throughout the simulation period. For large hydrogen plants, the scheduling of AWE systems, including on-standby-off state switching, should be further explored, along with a plant-wise performance comparison between the $N$-in-1 and 1-in-1 configurations. Incorporating uncertainty in energy input is another important direction.

\section*{Acknowledgement}

The authors gratefully acknowledge the financial support from the National Key Research and Development Program of China (2021YFB4000503) and the National Natural Science Foundation of China (52377116 and 52307126).

\section*{Declaration of Interest}

None.

\section*{Data Availability}

{\color{black}
The raw experimental dataset cannot be publicly released due to confidentiality restrictions from the plant owner and equipment manufacturer. The processed data are presented in graphical form in Appendix A. Access to de-identified or aggregated data may be considered with permission from the industrial partners.}
Other data related to this work are available upon request.

\break

\section*{Appendix A: Experimental Data of a 4,000 Nm$^3$/h-Rated 4-in-1 Electrolysis System}
\label{sec:exp}

\setcounter{figure}{0}
\renewcommand{\thefigure}{A\arabic{figure}}
\setcounter{table}{0}
\renewcommand{\thetable}{A\arabic{table}}

\begin{figure}[H]
	\centering
	\includegraphics[scale=0.94]{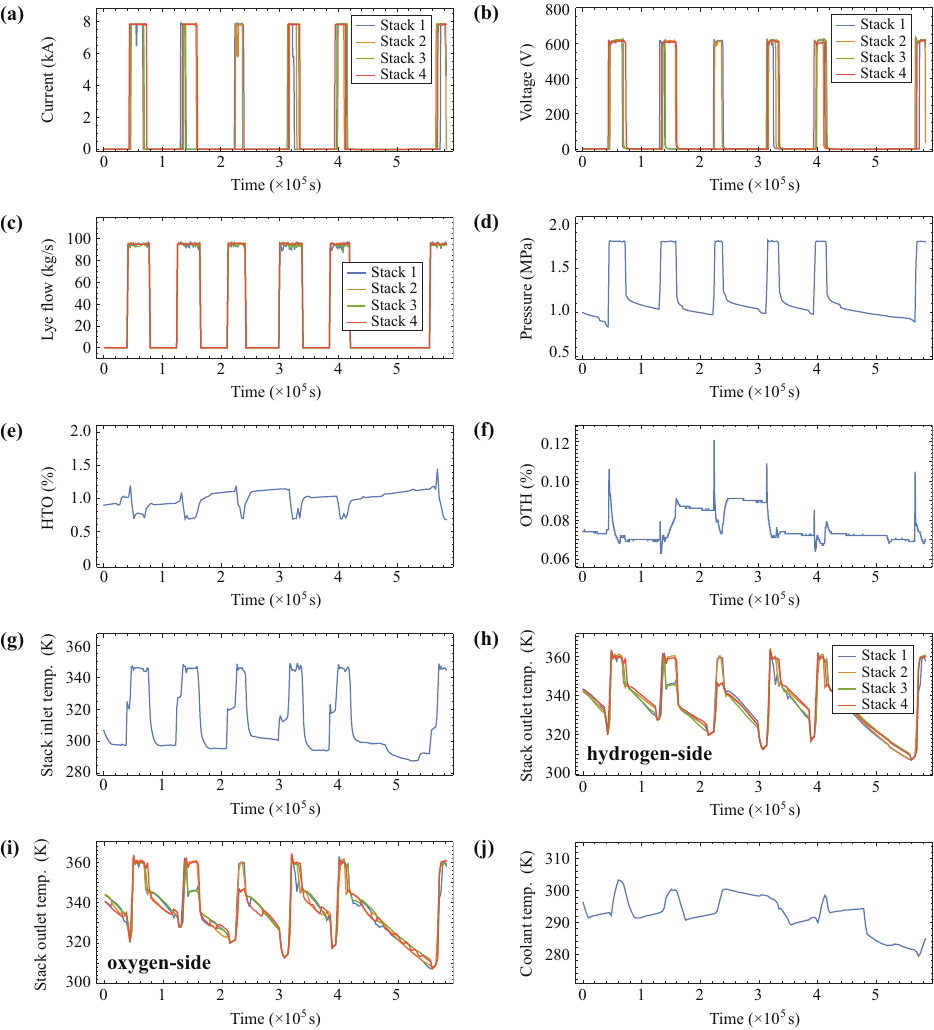}\vspace{-6pt}
	\caption{Experimental data of a 4,000 Nm$^3$/h-rated 4-in-1 alkaline electrolysis system used for model validation. (a) Stack current. (b) Stack voltage. (c) Lye flow rate of each stack. (d) System pressure. (e) HTO concentration in the gas phase of the lye-oxygen separator. {\color{black}(f) OTH concentration in the gas phase of the lye-hydrogen separator.} (g) Stack inlet temperature. {\color{black}(h) Hydrogen-side stack outlet temperatures. (i) Oxygen-side stack outlet temperatures.} (j) Coolant temperature.}
	\label{fig:expdata}
\end{figure}

\begin{figure}[H]
	\centering
	\includegraphics[scale=0.94]{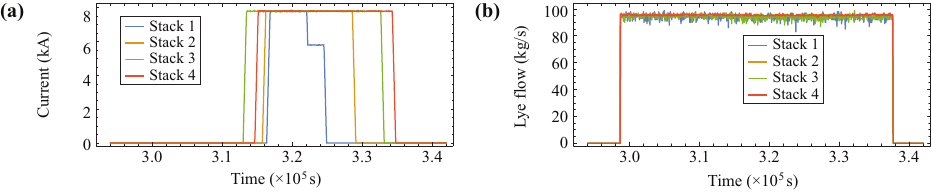}\vspace{-6pt}
	\caption{\color{black}Comparison of lye flow rate among 4 parallel stacks during a full startup-operation-shutdown cycle. (a) Electrolytic current. (b) Lye flow rate.}
	\label{fig:flowcomp}
\end{figure}

\begin{figure}[H]
	\centering
	\includegraphics[scale=0.94]{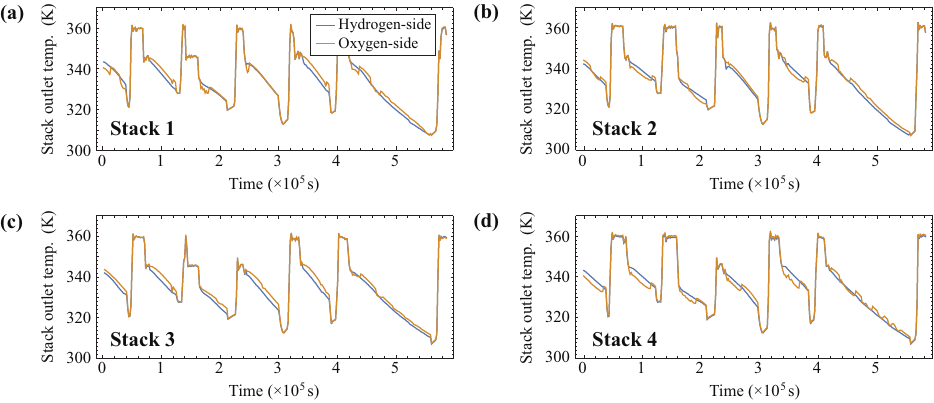}\vspace{-6pt}
	\caption{\color{black}Comparison between hydrogen-side and oxygen-side stack outlet temperatures. (a) Stack 1. (b) Stack 2. (c) Stack 3. (d) Stack 4. }
	\label{fig:tempcomp}
\end{figure}

The experimental data of a 4,000 Nm$^3$/h-rated 4-in-1 electrolysis system in a large-scale hydrogen plant in North China across 7 days are plotted in Fig. \ref{fig:expdata}. The data are used to validate the proposed state space model and calibrate model parameters in Section \ref{sec:validation}.

\section*{Appendix B: Manufacturer-Provided and Calibrated Parameters of the Alkaline Electrolysis System}
\label{sec:para}

\setcounter{figure}{0}
\renewcommand{\thefigure}{B\arabic{figure}}
\setcounter{table}{0}
\renewcommand{\thetable}{B\arabic{table}}

Table \ref{tab:parameter} lists the manufacturer-provided and calibrated parameters of the 4-in-1 AWE system used in the simulations in Section \ref{sec:case}.
{\color{black}
Some parameters, including geometric dimensions and rectifier ratings, are from factory datasheets.
Other parameters related to electrochemical, heat transfer, and mass transfer processes are calibrated using the Bayesian inference-based estimation method proposed in \cite{qiu2023dynamic}, which combines Markov chain Monte Carlo (MCMC) sampling and adaptive polynomial surrogate models, based on the 7-day experimental data presented in Section \ref{sec:experiment} and Appendix A.

Parameter calibration is conducted sequentially according to the physical hierarchy of the model. First, the electrochemical parameters are identified. Then, the thermal parameters are calibrated based on the measured temperature responses while fixing the electrochemical parameters. Finally, the HTO impurity dynamics parameters are estimated.
}

\begin{table}[H]\scriptsize
  \renewcommand{\arraystretch}{1.1}
  \caption{{\color{black}Parameters of the 4-in-1 alkaline water electrolysis system used in the simulation}}\vspace{6pt}
  \label{tab:parameter}
  \centering
  \begin{tabular}{ccc}
  \hline\hline
  Parameter                                                                     & & Value                  \\  \hline
  Number of stacks $N$                                                          & & 4 \\
  Rated hydrogen production                                                     & & $4,000$ Nm$^3$/h \\
  Rated hydrogen production per stack                                           & & $1,000$ Nm$^3$/h \\
  Rated electrolytic current per stack $I^\text{norm}$                          & & 7,800 A \\
  Maximal electrolytic current per stack $\overline{I}$                         & & 1.2$\times I^\text{norm}$ \\
  Maximal stack power $\overline{P}^{\text{ele}}$                               & & $6$  MW \\
  Number of cells per stack  $N^{\text{cell}}$                                  & & \color{black}$312$   \\
  Cell area $A^{\text{cell}}$                                                      & & $2$ m$^2$   \\
  System pressure $\rho$                                                        & & \color{black}1.8 MPa  \\
  Differential pressure between half-cells $\Delta \rho$                        & & $0.1\% \rho$ {\color{black}\cite{qi2021pressure}}  \\
  Electrochemical parameters $r_1$, $r_2$, $r_3$, $s$                           & & \color{black}$8.175\times10^{-6}$ $\Omega$, $2.136\times10^{-7}$ $\Omega$/K, $-8.656\times 10^{-12}$ $\Omega$/Pa, $7.024\times 10^{-2}$  \\
  Electrochemical parameters $t_1$, $t_2$, $t_3$                                & & \color{black}$-1.756\times10^{-1}$ $\Omega$, $79.25$ $\Omega\cdot$K, $32.56$ $\Omega\cdot$K$^2$ \\
  Faraday efficiency parameters $f_1$, $f_2$                                    & & \color{black}$50+ 2.5 T_s$ (A$^2$), $0.92-6.25\times10^{-6}T_s$   \\
  Upper limit of cell voltage $\overline{U}$                                    & & $2.1$ V \\
  Ramping limits $\overline{r}^{\text{H}_2,\text{prod}} / \underline{r}^{\text{H}_2,\text{prod} }$    & & $\pm 20$ Nm$^3$/(h$\cdot$s)         \\
  \color{black} Reversible voltage $U^{\text{rev}}$                                           & & \color{black}$1.229$ V \cite{ulleberg2003modeling}\\
  \color{black} Thermal neutral voltage $U^{\text{th}}$                                       & & \color{black}$1.481$ V \cite{ulleberg2003modeling} \\
  \hline
   \color{black}Dynamic viscosity of lye (30\%wt)              & &    \color{black}2.3$\times$10$^{-\text{3}}$ Pa$\cdot$s \cite{nistreference,hodges2023critical} \\
  \color{black}Dynamic viscosity of gaseous hydrogen			 & & 	\color{black}0.9$\times$10$^{-\text{5}}$ Pa$\cdot$s \cite{nistreference,hodges2023critical}\\
  \color{black}Dynamic viscosity of gaseous oxygen			 & &	\color{black}2.2$\times$10$^{-\text{5}}$ Pa$\cdot$s  \cite{nistreference,hodges2023critical} \\

  \hline
  Ambient temperature $T_{\text{am}}$                                           & & $298$ K ($25\ ^\circ$C) \\
  Coolant temperature $T_{\text{c,in}}$                                         & & $288$ K ($15\ ^\circ$C) \\
  Stack temperature reference $T^{\text{ref}}_{\text{s,out}}$                   & & $358$ K ($85\ ^\circ$C) \\
  Stack temperature limit $\overline{T}$                                        & & $363$ K ($90\ ^\circ$C) \\
  Stack heat capacity $C_{i,\text{s}}$                                          & & \color{black}$4.140\times10^7$ J/K    \\
  Separator heat capacity $C_{\text{sep}}$                                      & & \color{black}$2.922\times10^7$ J/K    \\
  Heat capacity of heat exchanger $C_{\text{he}}$                               & & \color{black}$2.224\times10^7$ J/K    \\
    \color{black} Stack diameter $\phi_{s}$  							& & \color{black}1.90 m \\
  \color{black} Stack heat dissipation area	$A_{\text{s,diss}}$		& & \color{black}50.44  m$^2$ \\				
  \color{black} Stack emissivity (nickel-plated) $\epsilon_\text{s}$ & & \color{black}0.25 \cite{touloukian1970thermal} \\
  \color{black} Stefan-Boltzmann constant   $\sigma$  				& & \color{black}$5.670\times 10^{-8}$ W/(m$^2\cdot$K$^4$) \cite{nistreference}\\
  Heat exchange area $A_{\text{c}}$                                             & & \color{black}$103.2$ m$^2$    \\
  \color{black} heat transfer coefficient $k$                                    & & \color{black} $980$ W/K$\cdot$m$^2$         \\
  Maximal cooling water flow rate $\overline{v}_\text{c}$                       & & 0.032 m$^3$/s        \\ \hline
  Anode-side half-cell lye volume $V^{\text{an}}_{i,\text{lye}}$                & & 2.5 m$^3$    \\
  Separator volume $V_{\text{sep}}$                                             & & 10.288 m$^3$    \\
  Limits of lye flow rate per stack  $\overline{v}_{\text{lye}}$, $\underline{v}_{\text{lye}}$  & & 0.0335 m$^3$/s, 0.0101 m$^3$/s    \\
  Thickness of the diaphragm $\delta$                                           & & 500$\times$10$^{-6}$ m \\
  Diffusion coefficient $D^{\rm{H_2}}_{\rm{eff}}$                               & & 8.569$\times$10$^{-10}$ m$^2$/s    \\
  Permeability coefficient $K^{\text{H}_2}_{\text{eff}}$                        & & 2$\times$10$^{-16}$  m$^2$  {\color{black}\cite{qi2021pressure}}  \\
  \color{black}Effective hydrogen carrying coefficient $S^{\text{H}_2}_{\text{eff}}$			& & \color{black} $1.540 \times 10^{-6}$ mol/(m$^3\cdot$Pa) \\
  Separation time constant $\tau_\text{sep}$                                    & & \color{black}$240$ s    \\
  \hline \hline
  \end{tabular}
\end{table}

{\color{black}
\section*{Appendix C: Constants and Matrices in the MIQP-based Controller}
\label{sec:para}

\setcounter{figure}{0}
\renewcommand{\thefigure}{C\arabic{figure}}
\setcounter{table}{0}
\renewcommand{\thetable}{C\arabic{table}}
\setcounter{equation}{0}
\renewcommand{\theequation}{C\arabic{equation}}

The polyhedral approximation matrix $\bm{A}$ and vector $\bm{b}$ for the stacks in the $4,000$ Nm$^3$/h system is given in (\ref{eq:mat}). Computational parameters used in the MIQP-based controller are summarized in Table \ref{tab:paracontrol}.
\begin{gather} \footnotesize
 \bm{A} = 10^{-3} \times \left( {\begin{array}{*{20}{c}}
  - 0.8257 &  - 0.0138 & 4.2644  \\
  - 0.6995 &  - 0.0678 & 2.6563  \\
  - 0.0980 &  \ \  2.8250 & 0.2724  \\
  - 0.1213 &  - 3.4951 & 0.5954  \\
  - 0.6101 &  - 0.5327 & 3.5723  \\
  - 0.8976 &  - 0.4399 & 4.0286  \\
  - 0.7899 &  - 1.6049 & 4.2766  \\
  - 0.8152 &  - 1.0185 & 3.3689  \\
  - 0.9217 &  - 1.8512 & 4.4626  \\
  - 0.8066 &  \ \  0.2918 & 3.3343  \\
  - 0.6156 &  - 2.3788 & 3.6415  \\
  - 0.8928 &  \ \  0.0317 & 4.2911
\end{array}} \right), \ \
    \bm{b}=
\left( {\begin{array}{*{20}{c}}
 - 0.6378 \\
 - 0.0141 \\
 - 1.7369 \\
 \ \  0.2377 \\
 - 0.8492 \\
 - 0.0454 \\
 - 0.2596 \\
 \ \  0.2731 \\
 \ \  0.2214 \\
 - 0.1737 \\
 - 0.2606 \\
 - 0.4030
\end{array}} \right). \label{eq:mat}
\end{gather}

\begin{table*}[h]\footnotesize
  \renewcommand{\arraystretch}{1.1}
  \caption{{\color{black}Parameters used in the MIQP-based controller}}\vspace{6pt}
  \label{tab:paracontrol}
  \centering
  \begin{tabular}{cc}
  \hline\hline
  Parameter                                                                     &  Value                  \\  \hline
  \color{black}Number of discretization bits $N^{\text{d}}$ for current                      &  \color{black}6 \\
  \color{black}Number of discretization bits $N^{\text{d}}$ for lye flow rate                &  \color{black}5 \\
  \color{black}Number of discretization bits $N^{\text{d}}$ for coolant flow rate            &  \color{black}6 \\
  \color{black}Big-M for current                                                             &  \color{black}12,000 (in A) \\
  \color{black}Big-M for temperature                                                         &  \color{black}443 (in K) \\
  \color{black}Big-M for molar quantity of HTO impurity                                      &  \color{black}520 (in mol)  \\
  \color{black}MIQP solver tolerance                                                         &  \color{black}$10^{-6}$  \\
  \color{black}MIP gap                                                                       &  \color{black}0.01  \\
  \hline\hline
  \end{tabular}
\end{table*}
}

\break


\end{document}